\RequirePackage{amsmath}
\RequirePackage{fix-cm}
\documentclass{svjour3}                     
\smartqed  
\usepackage{natbib}
\usepackage{xspace}
\usepackage{graphicx}
\usepackage{amssymb}
\usepackage{tabu}
\usepackage{cancel}
\usepackage{algorithm}
\usepackage{algorithmicx}
\usepackage{algpseudocode}
\usepackage{subcaption}
\usepackage{booktabs, mathtools}
\usepackage[dvipsnames]{xcolor}
\usepackage{url}
\usepackage[titletoc,title]{appendix}

\algdef{SE}[DOWHILE]{Do}{doWhile}{\algorithmicdo}[1]{\algorithmicwhile\ #1}%

\DeclareMathOperator{\Order}{{\mathcal O}}

\newcommand{\bm}[1]{\boldsymbol{#1}}
\newcommand{\mSigma}{\mathsf{\Sigma}}

\newcommand{\dif}[1]{\text{d}{#1}}
\newcommand{\D}[1]{\text{d}{#1}}
\newcommand{\trace}[1]{\textup{trace}{#1}}

\newcommand{\naturals}{\mathbb{N}}
\newcommand{\natzero}{\mathbb{N}_0}
\newcommand{\reals}{\mathbb{R}}
\newcommand{\integers}{\mathbb{Z}}

\newcommand{\Ex}{\mathbb{E}}

\newcommand{\rC}{\mathring{C}}

\newcommand{\rlambda}{\mathring{\lambda}}

\newcommand{\vDelta}{{\boldsymbol{\Delta}}}
\newcommand{\vlambda}{{\bm{\lambda}}}

\newcommand{\vtheta}{{\bm{\theta}}}
\newcommand{\vzeta}{{\bm{\zeta}}}

\newcommand{\va}{\bm{a}}
\newcommand{\vA}{\bm{A}}
\newcommand{\vb}{\bm{b}}
\newcommand{\vc}{\bm{c}}
\newcommand{\vC}{\bm{C}}

\newcommand{\vh}{\bm{h}}
\newcommand{\vf}{\bm{f}}

\newcommand{\vm}{\bm{m}}

\newcommand{\vt}{\bm{t}}
\newcommand{\vv}{\bm{v}}
\newcommand{\vV}{\bm{V}}
\newcommand{\vw}{\bm{w}}

\newcommand{\vx}{\bm{x}}

\newcommand{\dvx}{\dif{\bm{x}}}

\newcommand{\dvt}{\dif{\bm{t}}}
\newcommand{\vrho}{\bm{\rho}}

\newcommand{\vy}{\bm{y}}
\newcommand{\vY}{\bm{Y}}

\newcommand{\vz}{\bm{z}}

\newcommand{\dvz}{\dif{\bm{z}}}

\newcommand{\vPsi}{\boldsymbol{\Psi}}

\newcommand{\vone}{\bm{1}}

\newcommand{\mA}{\mathsf{A}}
\newcommand{\mC}{\mathsf{C}}
\newcommand{\mP}{\mathsf{P}}
\newcommand{\mK}{\mathsf{K}}
\newcommand{\rmC}{\mathring{\mathsf{C}}}
\newcommand{\mCtheta}{{\mathsf{C}_{\vtheta}}}
\newcommand{\mCInv}{\mathsf{C}^{-1}}
\newcommand{\cov}{{\textup{cov}}}
\newcommand{\var}{{\textup{var}}}
\newcommand{\opt}{{\textup{opt}}}

\newcommand{\mL}{\mathsf{L}}

\newcommand{\mLambda}{\mathsf{\Lambda}}

\newcommand{\mV}{\mathsf{V}}
\newcommand{\mW}{\mathsf{W}}

\newcommand{\calN}{\mathcal{N}}
\newcommand{\me}{\mathrm{e}}

\newcommand{\hmu}{\widehat{\mu}}
\newcommand{\hsigma}{\widehat{\sigma}}

\newcommand{\MLE}{\textup{EB}}
\newcommand{\full}{\textup{full}}
\newcommand{\GCV}{\textup{GCV}}
\newcommand{\errtol}{\varepsilon}

\newcommand{\err}{\textup{err}}
\newcommand{\code}[1]{\texttt{#1}}

\def\abs#1{\ensuremath{\left \lvert #1 \right \rvert}}
\newcommand{\norm}[2][{}]{\ensuremath{\left \lVert #2 \right \rVert}_{#1}}

\providecommand{\argmin}{\operatorname*{argmin}}

\newcommand\figref{Figure~\ref}
\newcommand\secref{Section~\ref}

\graphicspath{{.}{./figures/grey/}{D:/Mega/MyWriteupBackup/Sep_2ndweek_1/}}

\journalname{Automatic Bayesian Cubature}

\allowdisplaybreaks
\begin{document}
\setlength\abovedisplayskip{1pt}
\setlength{\belowdisplayskip}{1pt}

\title{Fast Automatic Bayesian Cubature Using Lattice Sampling
}
%



\author{R. Jagadeeswaran         \and
        Fred J. Hickernell 
}


\institute{R. Jagadeeswaran \at
              Department of Applied Mathematics, \\
              Illinois Institute of Technology \\
              10 W. 32nd St., Room 208,
              Chicago IL 60616\\
              \email{jrathin1@iit.edu}           
           \and
           Fred J. Hickernell \at
           Center for Interdisciplinary Scientific Computation and \\
           Department of Applied Mathematics \\
           Illinois Institute of Technology \\
           10 W. 32nd St., Room 208, 
           Chicago IL 60616
           \\
           \email{hickernell@iit.edu} 
}

\date{Received: date / Accepted: date}

\maketitle

\begin{abstract}
Automatic cubatures approximate integrals to user-specified error tolerances.  For high dimensional problems, it is difficult to adaptively change the sampling pattern, but one can automatically determine the
sample size, $n$, given a reasonable, fixed sampling pattern. We take this approach here using a Bayesian perspective.  We postulate that the integrand is an instance of a Gaussian stochastic process parameterized by a constant mean and a covariance kernel defined by a scale parameter times a parameterized function specifying how the integrand values at two different points in the domain are related.
These hyperparameters are inferred or integrated out using integrand values via one of three techniques:  empirical Bayes, full Bayes, or generalized cross-validation. The sample size, $n$, is increased until the half-width of the credible interval for the Bayesian posterior mean is no greater than the error tolerance. 

The process outlined above typically requires a computational cost of $O(N_{\text{opt}}n^3)$, where $N_{\text{opt}}$ is the number of optimization steps required to identify the hyperparameters. Our innovation is to pair low discrepancy nodes with matching covariance kernels to lower the  computational cost to $O(N_{\text{opt}} n \log n)$.   This approach is demonstrated explicitly with rank-1 lattice sequences and shift-invariant kernels.  Our algorithm is  implemented in the Guaranteed Automatic Integration Library (GAIL).

\end{abstract}

\section{Introduction}
\label{intro}
Cubature is the problem of inferring a numerical value for an integral, 
$ \mu := \int_{\reals^d} g(\vx) \, \dif \vx$, where $\mu$ has no closed form analytic expression. Typically, $g$ is accessible as a black-box algorithm. 
Cubature is a key component of many problems in scientific computing, finance, statistical modeling, and machine learning.  

The integral may often be expressed as
\begin{align}
\label{eqn:defn_mu}
\mu:= \Ex[f(\boldsymbol{X})] = \int_{[0,1]^d} f(\vx)\, \dif\vx, 
\end{align}
where $f:[0,1]^d \to \reals$ is the integrand, and $\boldsymbol{X} \sim \mathcal{U}[0,1]^d$.  The process of transforming the original integral into the form of \eqref{eqn:defn_mu} is not addressed here. See \cite[Section 2.11]{DicEtal14a} 
for a discussion of variable transformations. The cubature may be an affine function of integrand values:
\begin{align}
\label{eqn:defn_hmu}  
\hmu := w_0 + \sum_{i=1}^{n} f(\vx_i) w_i,
\end{align}
where the weights, $w_0$, and  $\vw = (w_i)_{i=1}^n \in \reals^n$, and the nodes, $\{\vx_i\}_{i=1}^n \subset [0,1]^d$, are chosen to make the error, $\abs{\mu - \hmu}$, small. The integration domain $[0,1]^d$ is convenient for the low discrepancy node sets \citep{DicEtal14a,SloJoe94} that we use.  The nodes are assumed to be deterministic.

Users of cubature algorithms typically want the error to be no greater than their specified error tolerance, denoted by $\varepsilon$.  That is, they want
\begin{align}
\label{eqn:err_crit} 
\abs{\mu - \hmu} \leq \errtol .
\end{align}
Some stopping criteria for choosing $n$ are heuristic.  Rigorous algorithms satisfying \eqref{eqn:err_crit}  typically require strong a priori assumptions about the integrand, such as an upper bound on its variance (for simple Monte Carlo) or total variation (for quasi-Monte Carlo).  We take a Bayesian approach by constructing a stopping criterion that is based on a credible interval.  We build upon the work of \cite{Dia88a}, \cite{OHa91a}, \cite{Rit00a}, \cite{RasGha03a}, \cite{BriEtal18a}, and others.  Our algorithm is an example of \emph{probabilistic numerics}.

Our primary contribution is to demonstrate how the choice of a family of covariance kernels that match the low discrepancy sampling nodes facilitates fast computation of the cubature and the data-driven stopping criterion.  Our Bayesian cubature requires a computational cost of
\begin{equation} \label{eqn:OuralgoCost}
    \Order\bigl(n \$(f) + N_{\opt}[n\$(C) + n \log(n)] \bigr),
\end{equation} 
where $\$(f)$ is the cost of one integrand value, $\$(C)$ is the cost of a single covariance kernel value,  $\Order(n \log(n))$ is the cost of a fast Fourier transform, and $N_{\opt}$ is an upper bound on the number of optimization steps required to choose the hyperparameters. If function evaluation is expensive, e.g., the output of a computationally intensive simulation, or if $\$(f) = \Order(d)$ for large $d$, then $\$(f)$ might be similar in magnitude to $N_{\opt} \log(n)$ in practice.  Typically, $\$(C) = \Order(d)$.  Note that the $\Order(n \log(n))$ contribution is $d$ independent.

By contrast to our fast algorithm, the typical computational cost for Bayesian cubature is
\begin{equation} \label{eqn:TheiralgoCost}
    \Order\bigl(n \$(f) + N_{\opt}[n^2\$(C) + n^3] \bigr),
\end{equation} 
which is explained in Section \ref{sec:bayes_cubature_algo}. Note that aside from evaluating the integrand, the computational cost in \eqref{eqn:TheiralgoCost} is much larger than that in \eqref{eqn:OuralgoCost}.  

\cite{Hic17a} compares different approaches to cubature error analysis depending on whether the rule is deterministic or random and whether the integrand is assumed to be deterministic or random.  Error analysis that assumes a deterministic integrand lying in a Banach space leads to an error bound that is typically impractical for deciding how large $n$ must be to satisfy \eqref{eqn:err_crit}.  The deterministic error bound includes a (semi-) norm of the integrand, often called the variation, which is often more complex to compute than the original integral.

\cite{HicJim16a} and \cite{JimHic16a} have developed stopping criteria for cubature rules based on low discrepancy nodes by tracking the decay of the discrete Fourier coefficients of the integrand.  The algorithm proposed here also relies on discrete Fourier coefficients, but in a different way.  Although we only explore automatic Bayesian cubature for absolute error tolerances, the recent work by \cite{HicEtal17a} suggests how one might accommodate more general error criteria, such as relative error tolerances.

Section \ref{sec:BC} explains the Bayesian approach to estimating the posterior cubature error and defines our automatic Bayesian cubature. Although much of this material is known, it is included for completeness.  We end Section \ref{sec:BC}  by demonstrating why Bayesian cubature is typically computationally expensive.
Section \ref{sec:fast_BC}  introduces the concept of covariance kernels that match the nodes, which expedites the computations required by our automatic Bayesian cubature. 
Section \ref{sec:shift_invariant_kernel} implements this concept for shift-invariant kernels and rank-1 lattice nodes.  This section also describes how to avoid cancellation error for covariance kernels of product form.  Numerical examples are provided in Section \ref{sec:NumExp} to demonstrate our new algorithm.  We conclude with a brief discussion.

\section{Bayesian Cubature} \label{sec:BC}

\subsection{Bayesian posterior cubature error}
\label{sec:BayesPostErr}

We assume that the integrand, $f$, is an instance of a Gaussian stochastic process  i.e., $f \sim \mathcal{GP}(m,s^2 C_\vtheta)$  \citep{Dia88a,OHa91a,Rit00a,RasGha03a,BriEtal18a}.  Specifically, $f$ is a real-valued random function with constant mean $m$ and covariance kernel $s^2C_\vtheta$:
\begin{gather*}
        m = \Ex[f(\vx)] \qquad \forall \vx \in \reals^d, \\
        \Ex\{[f(\vt) - m][f(\vx) -m]\} = s^2 C_\vtheta(\vt,\vx) \qquad \forall \vt, \vx \in \reals^d.
\end{gather*}
Here $s$ is a positive scale factor, and $C_\vtheta: [0,1]^d \times [0,1]^d \to \mathbb{R} $ is a symmetric, positive-definite function and parameterized by the vector $\vtheta$:
\begin{multline} \label{FJH:eq:CondPosDef}
\mC_\vtheta^T = \mC_\vtheta,  \; \va^T \mC_\vtheta \va > 0,  \text{ where }  \mC_\vtheta = \left(  C_\vtheta(\vx_i,\vx_j)  \right)_{i,j=1}^n,\\
 \forall \va \ne 0, \;
 n\in \mathbb{N}, \; \text{ distinct }\vx_1, \ldots, \vx_n \in [0,1]^d.
\end{multline}
Procedures for estimating or integrating out the hyperparameters $m$, $s$, and $\vtheta$ are explained later in this section.

Furthermore, for a Gaussian process, all vectors of linear functionals of $f$ have a multivariate Gaussian distribution. For any deterministic sampling scheme with distinct nodes, $\{\vx_i\}_{i=1}^n$, and defining  $\vf  := \left( f(\vx_i)\right)_{i=1}^n$ as the multivariate Gaussian vector of function values, it follows from the definition of a Gaussian process that 
\begin{subequations} \label{eqn:fGaussDist}
\begin{align}
\vf  & \sim \calN(m \vone, s^2 \mCtheta),  \qquad \text{where $\vone$ is a vector of all ones,}\\
\mu & \sim \calN(m, s^2 c_{0\vtheta}), 
\\
\text{where }
c_{0\vtheta} &:= \int_{[0,1]^{d}\times [0,1]^{d}} C_\vtheta(\vt,\vx) \, \dif{\vt} \, \dif{\vx}, \text{ and}\\
\cov(\vf, \mu) &= \left(  \int_{[0,1]^d} C_\vtheta(\vt,\vx_i) \, \D \vt \right)_{i=1}^n  =: \vc_\vtheta.
\end{align}
\end{subequations}
Here, $c_{0\vtheta}$ and $\vc_\vtheta$ depend explicitly on $\vtheta$.  We assume that $C_\vtheta$ is simple enough that the integrals in these definitions can be computed analytically.  We need the following lemma pertaining to a conditional Gaussian distribution to derive the distribution of the posterior error of our cubature. 

\begin{lemma} \citep[(A.6), (A.11--13)]{RasWil06a} \label{thrm:condDist} If $\vY = (\vY_1, \vY_2)^T \sim \calN (\vm,\mSigma)$, where $\vY_1$ and $\vY_2$ are random vectors of arbitrary length, and 
	\begin{gather*}
	\vm = \begin{pmatrix} \vm_1 \\ \vm_2 \end{pmatrix} = \begin{pmatrix} \Ex(\vY_1) \\ \Ex(\vY_2) \end{pmatrix}, \\
	\mSigma = \begin{pmatrix}
	\mSigma_{11} & \mSigma_{21}^T \\ 	\mSigma_{21} & \mSigma_{22}
	\end{pmatrix} =
	\begin{pmatrix}
	\var(\vY_{1}) & \cov(\vY_{1}, \vY_2) \\ 	\cov(\vY_2,\vY_{1}) & \var(\vY_{2}),
	\end{pmatrix} 
	\end{gather*}
	then 
	\begin{align*}
	\vY_1 \vert \vY_2 \; \sim \; \calN \bigl(\vm_1 + \mSigma_{21}^T \mSigma_{22}^{-1}(\vY_2 - \vm_2),  \; \mSigma_{11} - \mSigma_{21}^T \mSigma_{22}^{-1} \mSigma_{21} \bigr).
	\end{align*}
Moreover, the inverse of the matrix $\mSigma$ may be partitioned as
\begin{gather*}
\mSigma^{-1} = \begin{pmatrix} \mA_{11} & \mA_{21}^T \\ \mA_{21} & \mA_{22} \end{pmatrix}, \\
\mA_{11} = (\mSigma_{11} - \mSigma_{12} \mSigma_{22}^{-1} \mSigma_{21})^{-1}, \qquad 
\mA_{21} = -  \mSigma_{22}^{-1} \mSigma_{21} \mA_{11}, \\ 
\mA_{22} = \mSigma_{22}^{-1} + \mSigma_{22}^{-1} \mSigma_{21} \mA_{11} \mSigma_{21}^T \mSigma_{22}^{-1}.
\end{gather*}

\end{lemma}

It follows from Lemma \ref{thrm:condDist} that the \emph{conditional} distribution of the integral given observed function values, $\vf = \vy$ is also Gaussian:
\begin{align} \label{eqn:condInteg}
\mu | (\vf = \vy) \sim \calN \bigl(m (1 - \vc_\vtheta^T \mCtheta^{-1} \vone)  + \vc_\vtheta^T \mCtheta^{-1} \vy, 
\;
s^2(c_{0\vtheta}  -\vc_\vtheta ^T \mCtheta^{-1} \vc_\vtheta) \bigr).
\end{align}
The natural choice for  the cubature is the posterior mean of the integral, namely, 
\begin{align}
\label{eqn:BayesCub}
\widehat{\mu}  \vert ( \vf = \vy)
= m(1 - \vc_\vtheta^T  \mCtheta^{-1} \vone)
+ \vc_\vtheta^T \mCtheta^{-1} \vy,
\end{align}
which takes the form of \eqref{eqn:defn_hmu}.
Under this definition, the cubature error has zero mean and a variance depending on the choice of nodes:
\begin{align*}
(\mu-\hmu) | (\vf = \vy)
 \sim  \calN 
\left(
0, \;
s^2 (c_{0\vtheta} - \vc_\vtheta^T\mCtheta^{-1}\vc_\vtheta) 
\right).
\end{align*}
A credible interval for the integral is given by 
\begin{subequations} \label{eqn_prob_confidence_interval}
\begin{gather}
\mathbb{P}_f \left[
|\mu-\hmu| \leq \err_{\textup{CI}}
\right] = 99\%, \\
\err_{\textup{CI}} = 2.58 s \sqrt{c_{0\vtheta} - \vc_\vtheta^T\mCtheta^{-1}\vc_\vtheta}.
\end{gather}
\end{subequations}
Naturally, $2.58$ and $99\%$ can be replaced by other quantiles and credible levels.

\subsection{Hyperparameter estimation}
The credible interval in \eqref{eqn_prob_confidence_interval} suggests how our automatic Bayesian cubature proceeds.  Integrand data is accumulated until the width of the credible interval, $\err_{\textup{CI}}$, is no greater than the error tolerance.  As $n$ increases, one expects $c_{0\vtheta} - \vc_\vtheta^T\mCtheta^{-1}\vc_\vtheta$ to decrease for well-chosen nodes, $\{\vx_i\}_{i=1}^n$.

Note that $\err_{\textup{CI}}$ has no explicit dependence on the integrand values, even though one would intuitively expect that a larger integrand should imply a  larger $\err_{\textup{CI}}$.  This is because the hyperparameters, $m, s$, and $\vtheta$, have not yet been inferred from integrand data.  After inferring the hyperparameters, $\err_{\textup{CI}}$ does reflect the size of the integrand values. This section describes three approaches to hyperparameter estimation.

\begin{theorem} \label{thm:param} There are at least three approaches to estimating or integrating out the hyperparameters defining the Gaussian process from which the integrand is drawn: empirical Bayes, full Bayes, and generalized cross-validation.  Under these three approaches, we have the following:
\begin{align}
    \label{eqn_m_MLE}
m_\MLE &= \frac{\vone^T \mCInv_\vtheta \vy }{ \vone^T \mCInv_\vtheta \vone}, \qquad
m_{\textup{GCV}} = \frac{\vone^T \mC_\vtheta^{-2} \vy}{\vone^T \mC_\vtheta^{-2} \vone}, \\
\label{eqn_s2_MLE}
s^2_{\MLE} 
&= 
\frac{1}{n}
\vy^T 
\left[ \mCInv_\vtheta - 
\frac{ \mCInv_\vtheta \vone \vone^T \mCInv_\vtheta }{\vone^T\mCInv_\vtheta \vone}
\right] \vy, \\
\nonumber
\hsigma_{\full}^2 
& = \frac{1}{n-1}
\vy^T\left[ \mC_\vtheta^{-1} 
- \frac{ \mC_\vtheta^{-1} \vone\vone^T \mC_\vtheta^{-1}}{\vone^T \mC_\vtheta^{-1} \vone}  \right]\vy
\\ 
\label{eqn:sigma2_full}
& \qquad \qquad \times  \left[\frac{(1 - \vc_{\vtheta}^T \mC_\vtheta^{-1} \vone)^2}{\vone^T \mC_\vtheta^{-1} \vone} + (c_{0\vtheta}  -\vc_{\vtheta} ^T \mC_\vtheta^{-1} \vc_{\vtheta}) \right], \\
\nonumber
 s^2_{\textup{GCV}} & = \vy^T \left[\mC_\vtheta^{-2} - \frac{\mC_\vtheta^{-2} \vone \vone^T \mC_\vtheta^{-2}}{\vone^T \mC_\vtheta^{-2} \vone}  \right] \vy  \left[ \trace(\mC_\vtheta^{-1}) \right]^{-1}, \\
\label{eqn:thetaMLE}
\vtheta_\MLE
&= \argmin_{\vtheta} \biggl \{
\log\left(\vy^T 
\left[ \mC_\vtheta^{-1} - 
\frac{ \mCInv_\vtheta \vone \vone^T \mCInv_\vtheta }{\vone^T\mCInv_\vtheta \vone}
\right] \vy 
\right)  
 +  \frac{1}{n} \log(\det(\mC_\vtheta))
\biggr \}, \\
\label{vthetaGCV}
\vtheta_{\GCV} &= \argmin_\vtheta \biggl\{\log \left(  \vy^T \left[\mC^{-2}_\vtheta - \frac{\mC^{-2}_\vtheta \vone \vone^T \mC^{-2}_\vtheta}{\vone^T \mC^{-2}_\vtheta \vone}  \right] \vy \right)  
 - \log \left ( \trace(\mC^{-2}_\vtheta) \right ) \biggr\}, \\
\label{eqn:cubMLE}
\hmu_\MLE  &= \hmu_\full =
\left(
\frac{ (1 - \vone^T  \mCInv_\vtheta\vc_{\vtheta} )  \vone }{ \vone^T \mCInv_\vtheta \vone}   +  \vc_{\vtheta} 
\right)^T  \mCInv_\vtheta \vy, \\
\label{eqn:muCV}
\hmu_{\GCV}
& = \left(\frac{(1 - \vone^T  \mC_\vtheta^{-1}\vc_{\vtheta}) \mC_\vtheta^{-1} \vone}{\vone^T \mC_\vtheta^{-2} \vone} + \vc_{\vtheta} \right)^T \mC_\vtheta^{-1} \vy.
\intertext{The credible intervals widths, $\err_{\textup{CI}}$, are given by}
\label{eqn:err_MLEGCV}
\err_{\mathsf{x}} & = 2.58 s_{\mathsf{x}} \sqrt{c_{0\vtheta} - \vc_{\vtheta}^T\mC_\vtheta^{-1}\vc_{\vtheta} }, \qquad \mathsf{x} \in \{\MLE, \GCV\},  \\ 
\err_{\textup{full}} 
& = t_{n-1,0.995} \hsigma_{\textup{full}} > \err_\MLE. \label{FJH:eq:errFull}
\end{align}
The resulting credible intervals are then
\begin{align}
\label{eqn_prob_CI}
\mathbb{P}_f \left[
|\mu-\hmu_{\mathsf{x}}| \leq \err_{\mathsf{x}} \right]  = 99\%, \qquad \mathsf{x} \in \{\MLE, \full, \GCV\}.
\end{align}
Here $t_{n-1,0.995}$ denotes the $99.5$ percentile of a standard Student's $t$-distribution with $n-1$ degrees of freedom.  In the formulas above, $\vtheta$ is assumed to take on the values $\vtheta_{\MLE}$ or $\vtheta_{\GCV}$ as appropriate.
\end{theorem}

In the theorem above, note that if the original covariance kernel, $C_\vtheta$, is replaced by $b C_\vtheta$ for some positive constant $b$, the cubature, $\hmu$, the estimates of $\vtheta$, and the credible interval half-widths, $\err_{\mathsf{x}}$ for $\mathsf{x} \in \{\MLE, \full, \GCV\}$, all remain unchanged.  The estimates of $s^2$ are multiplied by $b^{-1}$, as would be expected. 

\subsubsection{Proof for Empirical Bayes}  \label{sec:MLE}
The empirical Bayes approach estimates the parameters, $m$, $s$, and $\vtheta$ via maximum likelihood estimation.  The log-likelihood function of the parameters given the integrand data $\vy$ is:
\begin{multline*}
l(s,m,\vtheta | \vy)
= -\frac{1}{2} s^{-2} (\vy-m\vone)^T\mCInv_\vtheta(\vy-m\vone) 
\\
 - \frac{1}{2} \log(\det(\mC_\vtheta)) - \frac{n}{2} \log(s^2) + \text{constants.}
\end{multline*}
Maximizing the log-likelihood first with respect to $m$ and then with respect to $s$ yields the values given in Theorem \ref{thm:param}.
To obtain $\vtheta_\MLE$, we substitute $m_\MLE$ and $s_\MLE$ into $l(s,m,\vtheta | \vy)$, which leads directly to the optimization problem in \eqref{eqn:thetaMLE}.

The empirical Bayes estimate of $\vtheta$ balances minimizing the covariance scale factor, $s^2_{\MLE}$, against minimizing  $\det(\mC_\vtheta)$. 
Under these estimates of the parameters, the cubature \eqref{eqn:BayesCub} and the credible interval \eqref{eqn_prob_confidence_interval} are explicitly written as in Theorem \ref{thm:param}.
The quantities $c_{0\vtheta}$, $\vc_\vtheta$, and $\mC_\vtheta$ are assumed implicitly to be based on $\vtheta = \vtheta_\MLE$.

\subsubsection{Proof for Full Bayes} \label{sec:fullBayes}
Rather than use maximum likelihood to determine $m$ and $s$,  one can treat them as hyperparameters with a non-informative, conjugate prior, namely $\vrho_{m,s^2}(\xi, \lambda) \propto 1/\lambda$. We want to compute $\rho_{\mu|\vf}(z | \vy)$, the conditional posterior density of $\mu$ given the data $\vf = \vy$.  This may be expressed as 
\begin{align*}
    \rho_{\mu|\vf}(z | \vy) 
= \int_{0}^\infty \int_{-\infty}^\infty 
\rho_{\mu | m, s^2, \vf}(z | \xi, \lambda , \vy)  
\rho_{m, s^2 | \vf}(\xi, \lambda | \vy)  \, \D \xi \D \lambda, 
\end{align*}
where $\rho_{m, s^2 | \vf}$ is the posterior density of the hyperparameters given the integrand data. Bayes Theorem tells us that $\rho_{m, s^2 | \vf} \propto \rho_{\vf | m, s^2} \, \vrho_{m,s^2}$, so 
\begin{align*}
\nonumber 
    \rho_{\mu|\vf}(z | \vy) 
& = \int_{0}^\infty \int_{-\infty}^\infty 
\rho_{\mu | m, s^2, \vf}(z | \xi, \lambda , \vy) 
\rho_{\vf | m, s^2} (\vy | \xi, \lambda) \, \vrho_{m,s^2}(\xi,\lambda)  \, \D \xi \D \lambda \\
\nonumber
& \propto \left( 1 +  \frac{(z - \hmu_{\MLE})^2}{(n-1) \hsigma_{\textup{full}}^2} \right)^{-n/2},
\end{align*}
where $\hsigma_{\textup{full}}^2$ is given in Theorem \ref{thm:param}, and the result above is derived in Appendix \ref{appendix:full_bayes}.

This means that $\mu \vert (\vf = \vy )$, properly centered and scaled, has a Student's $t$-distribution with $n-1$ degrees of freedom.   The estimated integral is the same as in the empirical Bayes case, $\hmu_{\textup{full}} = \hmu_{\MLE}$, but the credible interval is wider, as stated in the Theorem \ref{thm:param}.

Because the shape parameter, $\vtheta$, enters the definition of the covariance kernel in a non-trivial way, the only way to treat it as a hyperparameter and assign a tractable prior would be for the prior to be discrete.  We believe that choosing such a prior in practice involves too much guesswork, so we choose to use either $\vtheta_\MLE$ or $\vtheta_\GCV$.

\subsubsection{Proof for Generalized Cross-Validation} \label{sec:GCV}
A third parameter selection technique is \emph{leave-one-out cross-validation} (CV).  Let $\mathring{y}_i = \Ex[f(\vx_i ) | \vf_{-i} = \vy_{-i}]$, where the subscript $-i$ denotes the vector excluding the $i^{\text{th}}$ component.  This is the conditional expectation of $f(\vx_i )$ given the parameters $m$, $s$, and $\vtheta$, and all data but the function value at $\vx_i$.  The cross-validation criterion, which is to be minimized, is the sum of squares of the difference between these conditional expectations and the observed values:
\begin{align} \label{FJH:eq:CVA}
\textup{CV} = \sum_{i=1}^n (y_i - \mathring{y}_i)^2.
\end{align}

Let $\mA = \mC^{-1}_\vtheta$, let $\vzeta = \mA (\vy - m \vone)$, and partition $\mC_\vtheta$, $\mA$, and $\vzeta$ as
\begin{equation*}
\mC_\vtheta = \begin{pmatrix} c_{ii}  & \vC_{-i,i}^T \\  \vC_{-i,i} & \mC_{-i,-i}\end{pmatrix}, \qquad
\mA = \begin{pmatrix} a_{ii}  & \vA_{-i,i}^T \\  \vA_{-i,i} & \mA_{-i,-i}\end{pmatrix}, \qquad \vzeta = \begin{pmatrix} \zeta_i   \\  \vzeta_{-i} \end{pmatrix},
\end{equation*}
where the subscript $i$ denotes the $i^{\text{th}}$ row or column, and the subscript $-i$ denotes all rows or columns except the $i^{\text{th}}$. Following this notation, Lemma \ref{thrm:condDist} implies that 
\begin{align*}
\mathring{y}_i & = m + \vC^T_{-i,i} \mC_{-i,-i}^{-1} (\vy_{-i} -m \vone)  \\
\zeta_i  & = a_{ii}(y_i - m) + \vA_{-i,i}^T(\vy_{-i} - m \vone) \\
& = a_{ii}[(y_i - m) - \vC^T_{-i,i} \mC_{-i,-i}^{-1} (\vy_{-i} -m \vone)] \\
& = a_{ii}(y_i - \mathring{y}_i).
\end{align*}
Thus, \eqref{FJH:eq:CVA} may be re-written as 
\begin{align*} 
\textup{CV} = \sum_{i=1}^n \left(\frac{\zeta_i}{a_{ii}} \right)^2, \qquad \vzeta = \mC^{-1}_\vtheta(\vy - m \vone).
\end{align*}
The \emph{generalized cross-validation} criterion (GCV) replaces the $i^{\text{th}}$ diagonal element of $\mA$ in the denominator by the average diagonal element of $\mA$ \citep{CraWah79a,GolHeaWah79a,Wah90}:
\begin{equation*} 
\textup{GCV}
= \frac{\sum_{i=1}^n\zeta_i^2}{\left(\frac 1n \sum_{i=1}^n a_{ii} \right)^2} 
= \frac{(\vy - m\vone)^T \mC^{-2}_\vtheta (\vy - m \vone)}{\left(\frac 1n \trace(\mC_\vtheta^{-1}) \right)^2}.
\end{equation*}

The loss function $\textup{GCV}$ depends on $m$ and $\vtheta$, but not on $s$.  Minimizing the GCV  yields the formulae in Theorem \ref{thm:param} for $m_{\textup{GCV}}$ and $\vtheta_{\textup{GCV}}$.  
Plugging the value of $m_\GCV$ into \eqref{eqn:BayesCub} yields the formulae in Theorem \ref{thm:param} for $\widehat{\mu}_{\textup{GCV}}$.

An estimate for $s$ may be obtained by noting that by Lemma \ref{thrm:condDist},
\begin{align*}
\var[f(\vx_i ) | \vf_{-i} = \vy_{-i}] = s^2 a_{ii}^{-1}.
\end{align*}
Thus, we may estimate $s$ using an argument similar to that used in deriving the GCV and then substituting $m_{\textup{GCV}}$ for $m$:
\begin{align*}
s^2 &= \var[f(\vx_i ) | \vf_{-i} = \vy_{-i}] a_{ii} \\ 
& \approx \frac 1n \sum_{i=1}^n (y_i - \mathring{y}_i)^2a_{ii}
 = \frac 1n \sum_{i=1}^n \frac{\zeta_i^2}{a_{ii}} \\ 
 & \approx \frac{ \frac 1n \sum_{i=1}^n \zeta_i^2}{\frac 1n \sum_{i=1}^n a_{ii} } = \frac{(\vy - m\vone)^T \mC_\vtheta^{-2} (\vy - m \vone)}{ \trace(\mC_\vtheta^{-1}) } \\ 
 & \approx  \frac{(\vy - m_\GCV\vone)^T \mC_{\vtheta_\GCV}^{-2} (\vy - m_\GCV \vone)}{ \trace(\mC_{\vtheta_\GCV}^{-1}) } =: s^2_{\textup{GCV}}.
\end{align*}
After simplification, $s^2_{\textup{GCV}}$ defined above becomes the formula in Theorem \ref{thm:param}.

The credible interval based on GCV corresponds to \eqref{eqn_prob_confidence_interval} with the estimated $m$, $s$, and $\vtheta$.  This completes the proof of Theorem \ref{thm:param}.

\subsection{The automatic Bayesian cubature algorithm}
\label{sec:bayes_cubature_algo}
The previous section presents three credible intervals, \eqref{eqn_prob_CI}, for $\mu$, the desired integral.  Each credible interval is based on different assumptions about the hyperparameters $m$, $s$, and $\vtheta$.  We stress that one must estimate these hyperparameters or assume a prior distribution on them because the credible intervals are used as stopping criteria for our cubature rule.  Since a credible interval makes a statement about a typical function---not an outlier---one must try to ensure that the integrand is a typical draw from the assumed Gaussian stochastic process.

Our  Bayesian cubature algorithm increases the sample size until the width of the credible interval is small enough.  This is accomplished through successively doubling the sample size.  The steps are detailed in Algorithm \ref{algorithm1}.

We recognize that multiple applications of our credible intervals in one run of the algorithm is not strictly justified.  However, if our integrand comes from the middle of the sample space and not the extremes, we expect our automatic Bayesian cubature to approximate the integral within the desired error tolerance with high probability.  The example in the next subsection and the examples in Section \ref{sec:NumExp} support that expectation. We also believe that an important factor contributing to the occasional failure of our algorithm is unreasonable parameterizations of the stochastic process from which the integrand is hypothesized to be drawn.  Overcoming this latter challenge is a topic for future research.

\algnewcommand{\IIf}[1]{\State\algorithmicif\ #1\ \algorithmicthen\ }
\algnewcommand{\IElse}{\unskip\ \algorithmicelse\ }
\algnewcommand{\EndIIf}{\unskip\ \algorithmicend\ \algorithmicif}

\begin{algorithm}
\caption{Automatic Bayesian Cubature}\label{algorithm1}
  \begin{algorithmic}[1]
  	\Require a generator for the sequence
  	$\vx_1, \vx_2, \ldots$; 
  	a black-box function, $f$; 
  	an absolute error tolerance,
  	$\varepsilon>0$; the positive initial sample size, $n_0$;
  	the maximum sample size $n_{\textup{max}}$
  	
      \State $n \gets n_0, \; n' \gets 0, \; \err_{\textup{CI}} \gets \infty$
      
      \While{$\err_{\textup{CI}} > \varepsilon$ and $n \le n_{\textup{max}}$}
      
        \State Generate $\{ \vx_i\}_{i=n' + 1}^{n}$ and sample $\{f(\vx_i)\}_{i=n'+1}^{n}$
        \State Compute $\vtheta$ by \eqref{eqn:thetaMLE} or \eqref{vthetaGCV}
        \State Compute $\err_{\textup{CI}}$  according to \eqref{eqn:err_MLEGCV} or \eqref{FJH:eq:errFull}
        
       	\State	$n' \gets n, \; n \gets 2n'$
        
        \EndWhile
        
        \State Update sample size to compute $\hmu$, $n \gets n'$
        \State Compute $\hmu$, the approximate integral,   according to \eqref{eqn:cubMLE} or \eqref{eqn:muCV}
      \State \Return $\hmu$, $n$,  and $\err_{\textup{CI}}$
  \end{algorithmic}
\end{algorithm}

As described above, the computational cost of Algorithm \ref{algorithm1} is the sum of the following:
\begin{itemize}
	\item $\Order\bigl(n\$(f) \bigr)$ for the integrand data, where $\$(f)$ is the computational cost of a single $f(\vx)$; $\$(f)$ may be large if it is the result of an expensive simulation; $\$(f)$ is typically proportional to $d$;
	
	\item $\Order\bigl(N_{\opt} n^2 \$(C_\vtheta) \bigr)$ for the evaluation of the Gram matrix $\mC_{\vtheta}$, $N_\opt$ is the number of optimization steps required, and  $\$(C_\vtheta)$ is the computational cost of a single $C_\vtheta(\vt,\vx)$; $\$(C_\vtheta)$ is typically proportional to $d$; and
	
	\item $\Order\bigl(N_{\opt} n^3 \bigr)$ for the matrix inversions and determinant calculations; this cost is independent of $d$.
	
\end{itemize}
As we see in the example in the next section, this cost increases quickly as the $n$ required to meet the error tolerance increases.  This motivates the fast Bayesian cubature algorithm presented in Section \ref{sec:fast_BC}.

\subsection{Example with the Mat\'ern kernel} \label{MVN_example}

To demonstrate the automatic Bayesian cubature Algorithm \ref{algorithm1}, consider a Mat\'ern covariance kernel:
\begin{align*}
C_{\theta}(\vx, \vt) = \prod_{k=1}^d \exp(-\theta|\vx_k-\vt_k|)(1+\theta |\vx_k-\vt_k|),
\end{align*}
and Sobol' points as the nodes.  (Sobol' points are a typical space-filling design.)
Also, consider the integration problem of evaluating  \emph{multivariate Gaussian probabilities}:
\begin{align}
\label{eqn:GaussDef}
\mu = \int_{(\va,\vb)} \frac{\exp\bigl(- \frac 12 \vt^T \mSigma^{-1} \vt \bigr)}{\sqrt{(2 \pi)^{d'} \det(\mSigma)}} \, \dvt,
\end{align}
where $(\va,\vb)$ is a finite, semi-infinite or infinite box in $\reals^{d'}$.  This integral does not have an analytic expression for general $\mSigma$, so cubatures are required.  

\cite{Gen93} introduced a variable transformation to transform \eqref{eqn:GaussDef} into an integral on the unit cube.  Not only does this variable transformation accommodate domains that are (semi-)infinite, it also tends to smooth out the integrand better, which expedites the cubature.  Let $\mSigma= \mL \mL^T$ be the Cholesky decomposition where $\mL = (l_{jk})_{j,k=1}^d$ is a lower triangular matrix.  Iteratively define
\allowdisplaybreaks
\begin{gather*}
\alpha_1 = \Phi(a_1), \qquad \beta_1 = \Phi(b_1),
\\
\alpha_j(x_1,\dots,x_{j-1}) = 
\Phi
\left(
\frac{1}{l_{jj}} 
\left(
a_j - \sum_{k=1}^{j-1} l_{jk} \Phi^{-1}(\alpha_k + x_k(\beta_k-\alpha_k))
\right)
\right), 
\\
\hspace{8cm} j=2,\dots,d,
\\
\beta_j(x_1,\dots,x_{j-1}) = 
\Phi
\left(
\frac{1}{l_{jj}} 
\left(
b_j - \sum_{k=1}^{j-1} l_{jk} \Phi^{-1}(\alpha_k + x_k(\beta_k-\alpha_k))
\right)
\right), 
\\
\hspace{8cm} j=2,\dots,d,
\end{gather*}
\begin{equation}
\label{fGenzdef}
f_{\text{Genz}}(\vx) = \prod_{j=1}^d [\beta_j(\vx) - \alpha_j(\vx)],
\end{equation}
where $\Phi$ is the univariate cumulative standard Gaussian distribution function.  Then, $\mu = \int_{[0,1]^{d'-1}} f_{\text{Genz}}(\vx) \, \dvx$. This approach transforms a $d'$ dimensional integral into a $d=d'-1$ dimensional integral.

\begin{figure}
	\captionsetup[subfigure]{labelformat=empty}
	\centering
		\includegraphics[width=0.7\linewidth]{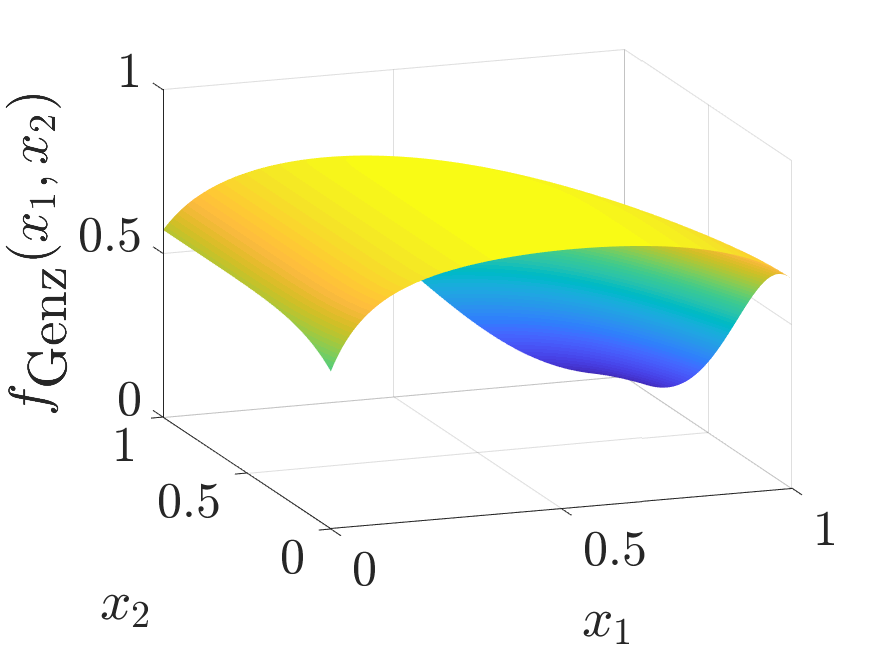}
	\caption{The $d'=3$ multivariate Gaussian probability transformed to an integral of $f_{\text{Genz}}$ of $d=2$. This plot can be reproduced using \code{IntegrandPlots.m} in GAIL.}
	\label{fig:MVN_Genz}
\end{figure}

\begin{figure}
	\centering
		\includegraphics[width=0.95\linewidth]{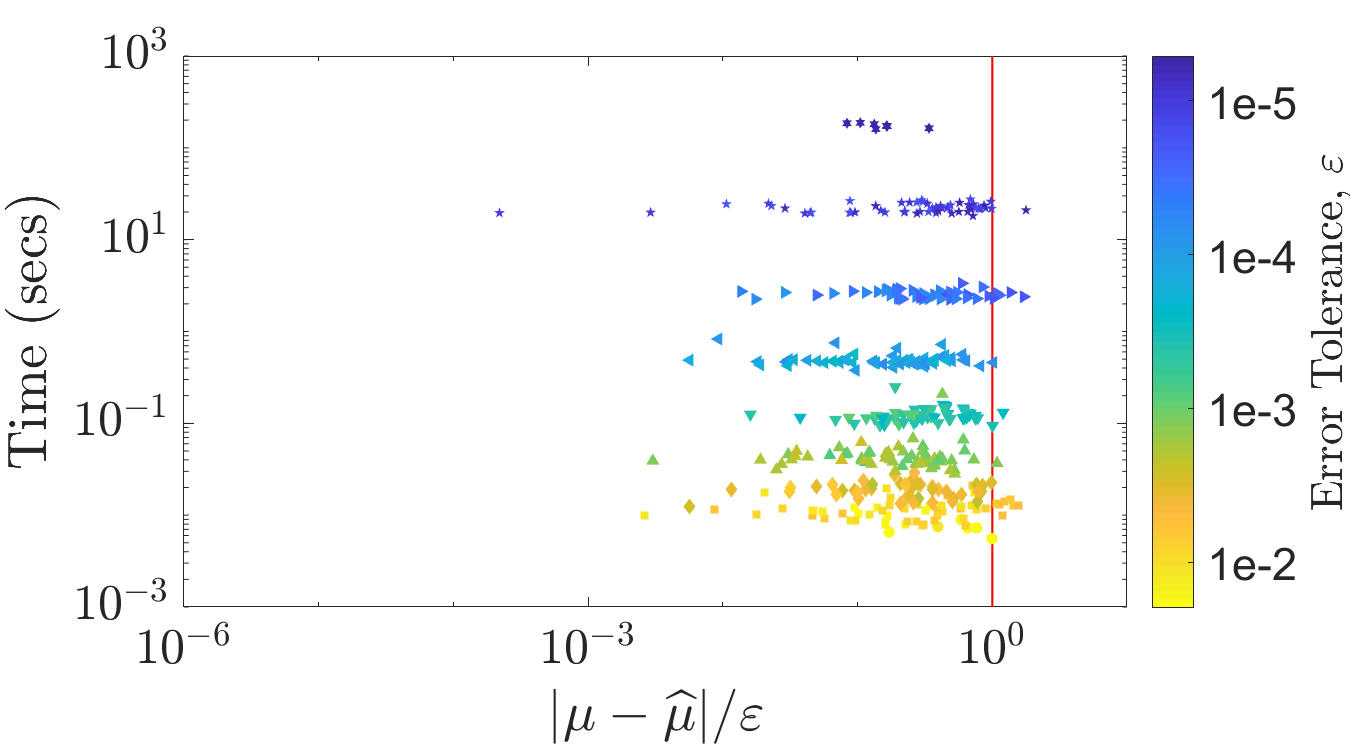}
	\centering
		\includegraphics[width=0.95\linewidth]{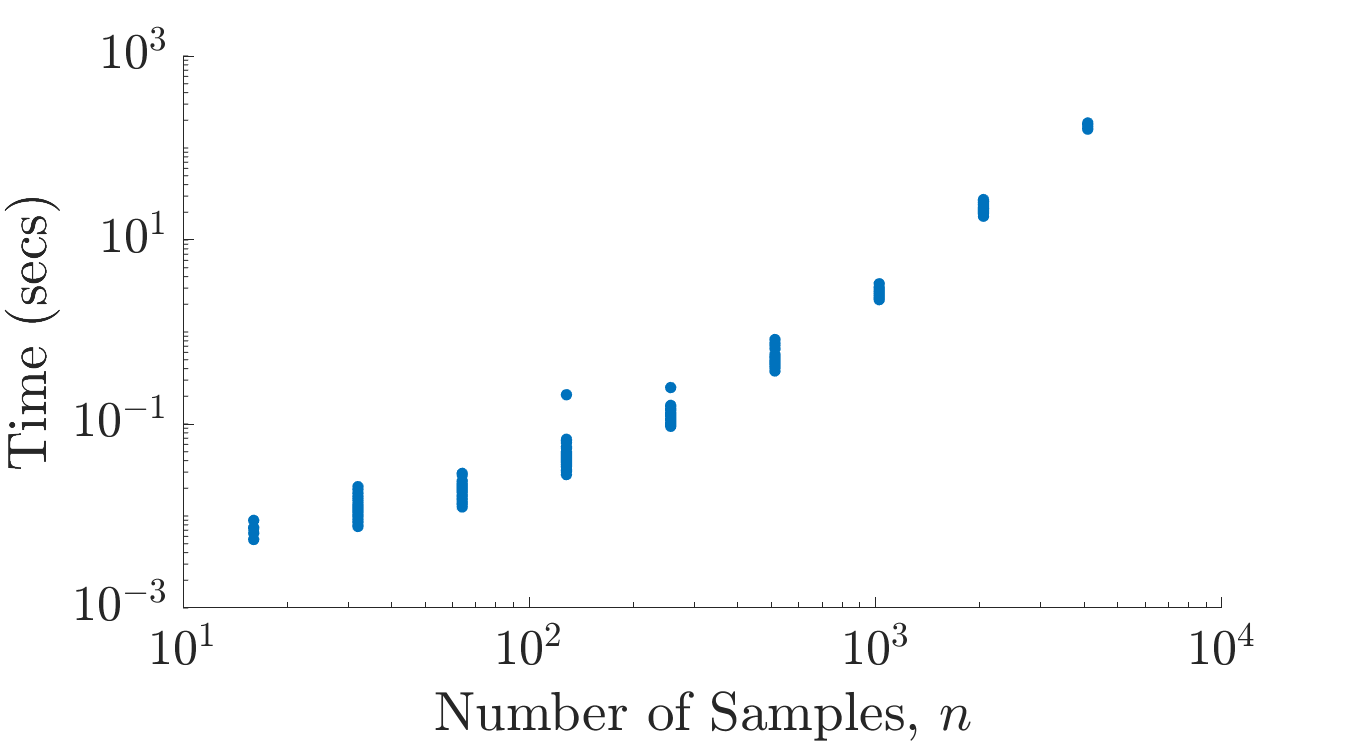}
	\caption{Multivariate Gaussian probability in $d=2$ estimated using Mat\'ern kernel and empirical Bayes stopping criterion. Top: Ratio of the integration error to the error-tolerance versus execution time. Bottom: Execution time rapidly increases with increasing $n$. These figures can be reproduced using \code{matern\_guaranteed\_plots.m} in GAIL.}
	\label{fig:MVN_Metern_d2b2}
\end{figure}

We use the following parameter values in the simulation: 
\begin{align*}
d' = 3, \quad \va = \begin{pmatrix}
-6 \\ -2 \\ -2
\end{pmatrix}, \quad 
\vb = \begin{pmatrix}
5 \\ 2 \\ 1
\end{pmatrix} , \quad 
\mL = \begin{pmatrix}
4 & 1 & 1 \\ 0 & 1 & 0.5 \\ 0 & 0 & 0.25
\end{pmatrix}.
\end{align*}
The node sets are randomly scrambled Sobol' points \citep{DicEtal14a,DicPil10a}.  The results are for 400 randomly chosen $\varepsilon$ in the interval $[10^{-5}, 10^{-2}]$ as shown in \figref{fig:MVN_Metern_d2b2}. In each run, the nodes are randomly scrambled.  The empirical Bayes credible intervals are used for stopping criteria.  We  observe that the algorithm meets the error criterion 95\% of the time even though we used 99\% credible intervals. One possible explanation is that the matrix inversions in the algorithm are ill-conditioned leading to numerical inaccuracies.  Another possible explanation is that this Mat\'ern covariance kernel is not a good match for the integrand.

As shown in \figref{fig:MVN_Metern_d2b2}, the computation time increases rapidly with $n$. 
The  empirical Bayes estimation of $\vtheta$, which requires repeated evaluation of the objective function, is the most time consuming of all. It takes tens of seconds to compute $\hmu_n$ with $\varepsilon = 10^{-5}$.   In contrast, this example in Section \ref{sec:NumExp} take less than a hundredth of a second to compute $\hmu_n$ with the same $\varepsilon$ using our new algorithm. Not only is the Bayesian cubature with the Mat\'ern kernel slow, but also $\mC_\vtheta$ becomes highly ill-conditioned as $n$ increases.
So, Algorithm \ref{algorithm1}  in its current form is impractical when $n$ must be large. 

\section{Fast Automatic Bayesian Cubature}\label{sec:fast_BC}

The generic automatic Bayesian cubature algorithm described in the previous section requires $\Order\bigl(n \$(f) +  N_{\opt} [n^2 \$(C_\vtheta) + n^3]\bigr)$ operations to compute the cubature. Now we explain how to speed up the calculations. A key is to choose covariance kernels that match the nodes, $\{\vx_i\}_{i=1}^n$, so that the vector-matrix operations required by Bayesian cubature can be accomplished using fast Bayesian transforms at a computational cost of $\Order\bigl(n \$(f) + N_{\opt} [n \$(C_\vtheta)  + n \log(n)] \bigr)$.

\subsection{Fast Bayesian Transform Kernel}
We make some assumptions about the relationship between the covariance kernel and the nodes.  In Section \ref{sec:shift_invariant_kernel} these assumptions are shown to hold  for rank-1 lattices and shift-invariant kernels.  Although the integrands and covariance kernels are real, it is convenient to allow related vectors and matrices to be complex.  A relevant example is the fast Fourier transform (FFT) of a real-valued vector, which is a complex-valued vector. 

We introduce some further notation:
\begin{align}
\mC = \mC_\vtheta &= \Big(C_\vtheta(\vx_i,\vx_j)\Big)_{i,j=1}^n  = (\vC_1,\dots,\vC_n) 
\label{eqn:ftk_factor}
= \frac 1n \mV \mLambda \mV^H , 
\quad \quad \mV^H = n \mV^{-1}, \\
\nonumber
\mV &= (\vv_1,\ldots,\vv_n)^T = (\vV_1,\ldots,\vV_n),
\qquad
\mC^p  = \frac 1n \mV \mLambda^{p} \mV^H, \qquad \forall p \in \integers.
\end{align}
In this and later sections, we drop the $\vtheta$ dependence of various quantities for simplicity of notation.  Here, $\mV^H$ is the Hermitian of $\mV$, $\vC_1,\dots,\vC_n$ are columns of $\mC$,  $\vV_1,\dots,\vV_n$ are columns of $\mV$, and $\vv_1,\dots,\vv_n$ are rows of $\mV$.  The normalization of $\mV$ assumed in \eqref{eqn:ftk_factor} conveniently allows the first eigenvector, $\vV_1$, to be the vector of ones in \eqref{fastcompAssumpB} below.  The columns of $\mV$ are eigenvectors of $\mC$, and $\mLambda$ is a diagonal matrix of eigenvalues of $\mC$.
For any $n \times 1$ vector $\vb$, define the notation  $\widetilde{\vb} := \mV^H \vb$.

We make three assumptions that facilitate  fast computation:
\begin{subequations} \label{fastcompAssump}
	\begin{gather}
	\label{fastcompAssumpA}
	\mV \text{ may be identified analytically}, \\
	\label{fastcompAssumpB}
	\vv_1 = \vV_1 = \vone, \\
	\label{fastcompAssumpC}
	\mV^H \vb  \text{ requires only $\Order(n \log(n))$ operations } \forall \vb.
	\end{gather}
\end{subequations}
We call the transformation $\vb \mapsto \mV^H \vb$ a \emph{fast Bayesian transform} and $C_\vtheta$ a \emph{fast Bayesian transform kernel}.  

Under assumptions \eqref{fastcompAssump} the eigenvalues may be identified as the fast Bayesian transform of the first column of $\mC$:
\begin{equation}
\vlambda 
= \begin{pmatrix}
\lambda_1 \\ \vdots \\ \lambda_n
\end{pmatrix} = \mLambda \vone = \mLambda \vv_1^* 
= \underbrace{\left( \frac 1n \mV^H  \mV \right) }_{\mathsf{I}} \mLambda \vv_1^* 
= \mV^H \left( \frac 1n \mV \mLambda \vv_1^* \right)
= \mV^H \vC_1 =  \widetilde{\vC}_1,
\label{eqn:fast_transform_to_eigvalues}
\end{equation}
where $\mathsf{I}$ is the identity matrix and $\vv_1^*$ is the complex conjugate of the first row of $\mV$.
Also note that the fast Bayesian transform of $\vone$ has a simple form
\begin{align*} 
\widetilde{\vone}
& = \mV^H \vone = \mV^H \vV_1 = 
\left(n, 0, \dots, 0 \right)^T.
\label{eqn:fast_transform_one}
\end{align*}

Many of the terms that arise in the calculations in  Algorithm \ref{algorithm1} take the form $\va^T\mC^{p}\vb$ for real $\va$ and $\vb$, and integer $p$.  These can be calculated via the transforms $\widetilde{\va} = \mV^H \va$ and $\widetilde{\vb} = \mV^H \vb$ as 
\begin{equation*}
\va^T\mC^p\vb = \frac 1n \va^T \mV \mLambda^p \mV^H \vb
= \frac 1n \widetilde{\va}^H\mLambda^p \widetilde{\vb}
= \frac 1n \sum_{i=1}^n \lambda_i^p \widetilde{a}_i^* \widetilde{b}_i.
\end{equation*}
Note that $\widetilde{\va}^*$ appears on the right side of this equation because $\va^T \mV = (\mV^H \va)^* = \widetilde{\va}^*$.  In particular,
\begin{align*}
\vone^T\mC^{-p}\vone & = \frac{n}{\lambda_1^p},
&
\vone^T\mC^{-p}\vy &= \frac{\widetilde{y}_1}{\lambda_1^p},
\\
\vy^T\mC^{-p} \vy &= \frac 1n \sum_{i=1}^n \frac{\abs{\widetilde{y}_i}^2}{\lambda_i^p},
&
\vc^T\mCInv \vone &= \frac{\widetilde{c}_1}{\lambda_1},\\
\vc^T\mCInv \vy &= \frac 1n \sum_{i=1}^n \frac{\widetilde{c}_i^* \widetilde{y}_i}{\lambda_i}, & 
\vc^T\mCInv \vc &= \frac 1n \sum_{i=1}^n \frac{\abs{\widetilde{c}_i}^2}{\lambda_i},
\end{align*}
where $\widetilde{\vy} = \mV^H \vy$ and 
$\widetilde{\vc} = \mV^H \vc$.  For any real $\vb$, with $\widetilde{\vb} = \mV^H\vb$, it follows that $\widetilde{b}_1$ is real since the first row of $\mV^H$ is $\vone$.

The covariance kernel used in practice also may satisfy an additional assumption:
\begin{align} \label{addAssump}
\int_{[0,1]^d} C_{\vtheta}(\vt,\vx) \, \D \vt = 1 \qquad \forall \vx \in [0,1]^d,
\end{align}
which implies that $c_{0\vtheta} = 1$ and $\vc_\vtheta = \vone$.  Under \eqref{addAssump}, the expressions above may be further simplified:
\begin{align*}
\vc^T\mCInv \vone =
\vc^T\mCInv \vc = \frac{n}{\lambda_1}.
\end{align*}
The assumptions and derivations above lead to the following theorem.

\begin{theorem} \label{thm:fastparam}
Under assumptions \eqref{fastcompAssump}, the parameters and credible interval half-widths in Theorem \ref{thm:param} may be expressed in terms of the fast Bayesian transforms of the integrand data, the first column of the Gram matrix, and $\vc_\vtheta$ as follows:
\begin{align}
\nonumber
m_\MLE &=  m_{\full} = m_{\GCV} =  \frac{\widetilde{y}_1}{n} = \frac 1n \sum_{i=1}^n y_i,
\\
\nonumber
s^2_\MLE 
& =
\frac{1}{n^2} 
\sum_{i=2}^n \frac{\abs{\widetilde{y}_i}^2}{\lambda_i}, \\
\nonumber
\widehat{\sigma}^2_{\textup{full}}
& =
\frac{1}{n(n-1)} \sum_{i=2}^n \frac{\abs{\widetilde{y}_i}^2}{\lambda_i}
\times
\left[\frac{\lambda_1}{n}{\left(1 - \frac{\widetilde{c}_1}{\lambda_1}\right)^2} + \left(c_0  - \frac 1n \sum_{i=1}^n \frac{\abs{\widetilde{c}_i}^2}{\lambda_i}\right) \right], \\
\nonumber 
s^2_{\textup{GCV}} & =  \frac 1{n} \sum_{i=2}^n \frac{\abs{\widetilde{y}_i}^2}{\lambda_i^2}  \left [ \sum_{i=1}^n \frac{1}{\lambda_i} \right]^{-1},
\end{align}
\begin{subequations}
\label{fastTheta}
\begin{align}
\label{eqn_MLE_loss_func_optimized_2}
\vtheta_\MLE
&= 
\argmin_{\vtheta}
\left[
\log\left(
\sum_{i=2}^n \frac{\abs{\widetilde{y}_i}^2}{\lambda_{\vtheta, i}}
\right) 
 + \frac{1}{n}\sum_{i=1}^n \log(\lambda_{\vtheta,i})
\right],\\
\label{thetaGCV} 
\vtheta_{\GCV} 
&= \argmin_\vtheta \left[ \log \left ( \sum_{i=2}^n \frac{\abs{\widetilde{y}_i}^2}{\lambda_{\vtheta, i}^2} 
\right)  
-2\log\left( \sum_{i=1}^n \frac{1}{\lambda_{\vtheta, i}} \right)
\right], 
\end{align}
\end{subequations}
\begin{align}
\nonumber
\hmu_\MLE  &= \hmu_{\full} = \hmu_{\GCV} =
\frac{\widetilde{y}_1}{n} +
\frac 1n \sum_{i=2}^n \frac{ \widetilde{c}_i^* \widetilde{y}_i}{\lambda_i}, \\
\nonumber
\err_\MLE  &
=
\frac{2.58}{n}\sqrt{
	\sum_{i=2}^{n} \frac{\abs{\widetilde{y}_i}^2}{\lambda_i}  
	\,
	\left( c_0 - \frac 1n \sum_{i=1}^n \frac{\abs{\widetilde{c}_i}^2}{\lambda_i} \right) 
}, \\
\nonumber
\err_{\full} & = t_{n-1,0.995} \hsigma_{\textup{full}}, \\
\nonumber
\err_{\textup{GCV}} & =
\frac{2.58}{n}\left\{\sum_{i=2}^n \frac{\abs{\widetilde{y}_i}^2}{\lambda_i^2}  \left [ \frac 1n \sum_{i=1}^n \frac{1}{\lambda_i} \right]^{-1} 
\left( c_0 - \frac 1n \sum_{i=1}^n \frac{\abs{\widetilde{c}_i}^2}{\lambda_i} \right) 
\right\}^{1/2}.
\end{align}
Under the further assumption \eqref{addAssump}, it follows that 
\begin{equation}
\label{muhatGCV-FB-MLE-Simple}
\hmu_\MLE  = \hmu_{\full} = \hmu_{\GCV} =
\frac{\widetilde{y}_1}{n} = \frac 1n \sum_{i=1}^n y_i,
\end{equation}
and so $\hmu$ is simply the sample mean.  Also, under assumption \eqref{addAssump}, the credible interval half-widths simplify to
\begin{subequations}
\label{errSimple}
\begin{align}
\label{eq:errMLEAllAsump}
\err_\MLE  &
=
\frac{2.58}{n}\sqrt{
	\sum_{i=2}^{n} \frac{\abs{\widetilde{y}_i}^2}{\lambda_i}  
	\,
	\left( 1 -  \frac{n}{\lambda_1} \right) 
}, \\
\label{FJH:eq:errFullSimple}
\err_{\textup{full}}
&=
t_{n-1,0.995}
\sqrt{\frac{1}{n(n-1)} \sum_{i=2}^n \frac{\abs{\widetilde{y}_i}^2}{\lambda_i}  \left(\frac{\lambda_1}{n}  - 1  \right)}, \\
\nonumber
\err_{\textup{GCV}} & =
\frac{2.58}{n}\left\{\sum_{i=2}^n \frac{\abs{\widetilde{y}_i}^2}{\lambda_i^2}  \left [ \frac 1n \sum_{i=1}^n \frac{1}{\lambda_i} \right]^{-1} 
\left( 1 -  \frac{n}{\lambda_1} \right)  
\right\}^{1/2}. \label{errGCVSimple}
\end{align}
\end{subequations}
In the formulas for the credible interval half-widths, $\vlambda$ depends on $\vtheta$, and  $\vtheta$ is assumed to take on the values $\vtheta_{\MLE}$ or $\vtheta_{\GCV}$ as appropriate.
\end{theorem}

\section{Integration Lattices and Shift Invariant Kernels}
\label{sec:shift_invariant_kernel}

The preceding sections lay out an automatic Bayesian cubature algorithm whose computational cost is drastically reduced.  However, this algorithm relies on covariance kernel functions, $C_\vtheta$ and node sets, $\{\vx_i\}_{i=1}^n$ that satisfy assumptions \eqref{fastcompAssump}.  We  also want to satisfy assumption \eqref{addAssump}.  
To conveniently facilitate the fast Bayesian transform, it is assumed in this section and the next that $n$ is power of $2$.  

\subsection{Extensible Integration Lattice Node Sets}

We choose a set of nodes defined by a shifted extensible integration lattice node sequence, which takes the form
\begin{equation} \label{eqn:lattice_def}
\vx_{i} = \vh \phi(i-1) + \vDelta \mod \vone, \qquad i \in \naturals.
\end{equation} 
Here, $\vh$ is a $d$-dimensional generating vector of positive integers, $\vDelta$ is some point in $[0,1)^d$, often chosen at random, and $\{\phi(i)\}_{i=0}^\infty$ is the van der Corput sequence, defined by reflecting the binary digits of the integer about the decimal point, i.e., 
\begin{align} \label{vdCDef}
\begin{array}{r|ccccccccccccc}
i & 0 & 1 & 2 & 3 & 4 &  5 & 6 & 7 & \cdots \\
i & 0_2 & 1_2 & 10_2 & 11_2 & 100_2 & 101_2 & 110_2 & 111_2  & \cdots\\
\toprule
\phi(i) & {}_2.0 &  {}_2.1 & {}_2.01 &  {}_2.11  & {}_2.001 &  {}_2.101 & {}_2.011 &  {}_2.111 & \cdots\\
\phi(i) & 0 &  0.5 &  0.25 & 0.75 &  0.125 & 0.625  &  0.375 & 0.875 & \cdots
\end{array}
\end{align}
Note that 
\begin{align} \label{phiprop}
n\phi:\{0, \ldots, n-1 \} \to \{0, \ldots, n-1\} \quad
\text{is one-to-one},
\end{align}
assuming $n$ is a power of $2$.

\begin{figure}[htp]
	\centering
	\includegraphics[height=5cm]{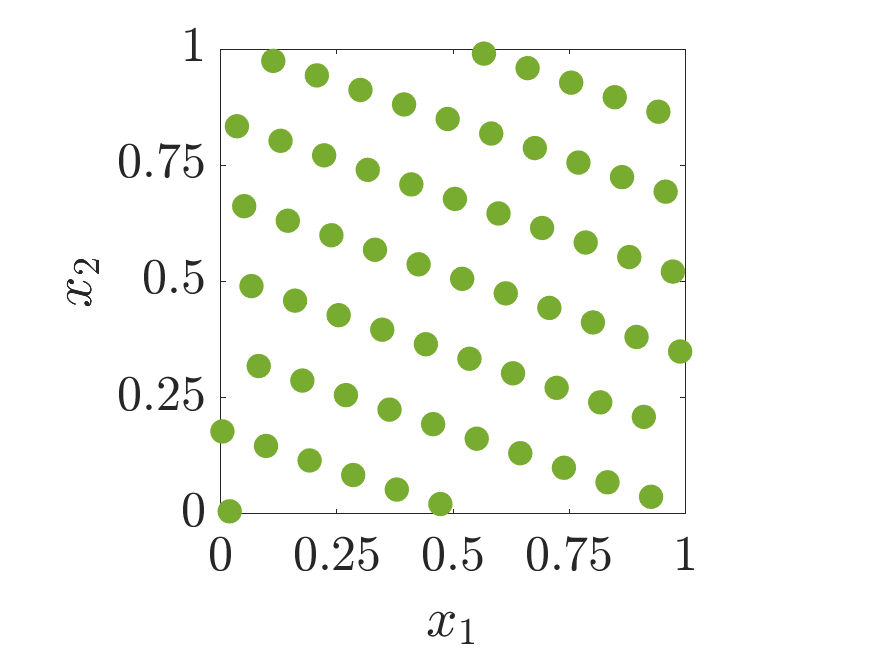}
	\caption{Example of a shifted integration lattice node set  in $d=2$. 
	This figure can be reproduced using \code{PlotPoints.m} in GAIL} \label{latticefig}
\end{figure}

An example of $64$ nodes is given in \figref{latticefig}.  The even coverage of the unit cube is ensured by a well chosen generating vector, $\vh$.  The choice of generating vector is typically done offline by computer search.  See \cite{HicNie03a}  and \cite{DicEtal14a} for more on extensible integration lattices.

\subsection{Shift Invariant Kernels}
The covariance kernels $C_\vtheta$ that match integration lattice node sets have the form
\begin{subequations} \label{eqn:ourkernel}
\begin{align} \label{eq:shInv}
C_\vtheta(\vt,\vx) = K_\vtheta(\vt - \vx \bmod \vone).
\end{align}
This is called a \emph{shift-invariant kernel} because shifting both arguments of the covariance kernel by the same amount leaves the value unchanged.  Here, $K_\vtheta$ is periodic and must  be of the form that ensures that $C_\vtheta$ is symmetric and positive definite, as assumed in \eqref{FJH:eq:CondPosDef}. 

A family of shift-invariant kernels is constructed via even degree Bernoulli polynomials:
\begin{multline}
\label{the_kernel_eqn_bernoulli}
K_\vtheta(\vx) =
\prod_{l=1}^d \biggl[
1 - (-1)^{r} \eta B_{2r}( x_l) \biggr], \\  
\forall \vt,\vx \in [0,1]^d, \  \vtheta = (r,\eta), \ r \in \naturals, \ \eta > 0.
\end{multline}
\end{subequations}
Symmetric, periodic, positive definite kernels of this form appear in \cite{Hic96a} and  \cite{DicEtal14a}.  Bernoulli polynomials are described in \cite[Chapter 24]{OlvEtal10a}.

\begin{figure}
	\centering
	\includegraphics[width=0.9\linewidth]{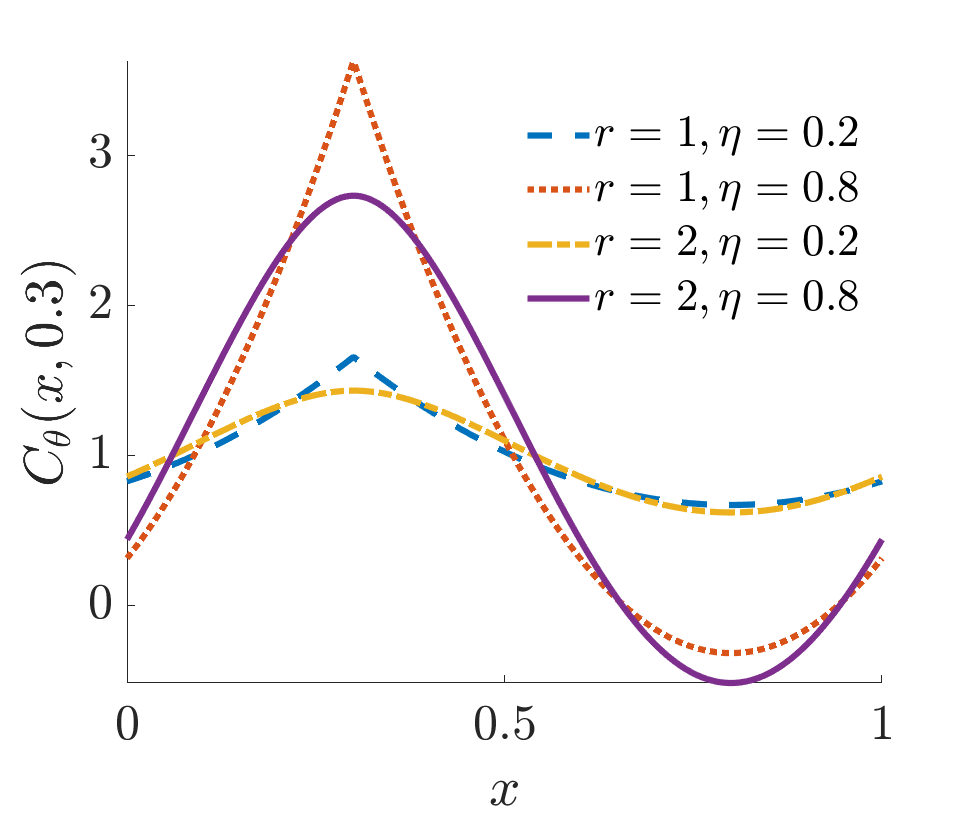}
	\caption[Fourier kernel]{Shift invariant kernel in 1D shifted by 0.3 to show the discontinuity. This figure can be reproduced using \code{plot\_fourier\_kernel.m} in GAIL. }
	\label{fig:fourierkernel-dim1}
\end{figure}

Larger $r$ implies a greater degree of smoothness of the covariance kernel.  Larger $\eta$ implies greater fluctuations of the output with respect to the input.  Plots of $C_\vtheta(\cdot, 0.3)$ are given in \figref{fig:fourierkernel-dim1} for various $\vtheta=(r, \eta)$ values.

Lattice cubature rules are known to have convergence rates that depend on the smoothness of the integrands, but that are rather independent of the choice of the integration lattice \citep{DicEtal14a}.  Thus, we expect integration lattice node sets to perform well regardless of the smoothness of the covariance kernel.  The bigger concern is whether the derivatives of the integrand are as smooth as the covariance kernel implies.  This topic is touched upon again in Section  \ref{period_var_tx}.

\subsection{The Gram Matrix as the Permutation of a Circulant Matrix}
For general shift-invariant covariance kernels, the Gram matrix takes the form of a permutation of the rows and columns of a circulant matrix. By the properties of $\phi$ in \eqref{phiprop}, it follows that
\begin{equation} \label{PermMat}
    \mP = \bigl( \delta_{n\phi(i-1), j-1}  \bigr)_{i,j=1}^n
\end{equation}
is a permutation matrix, where $\delta_{\cdot,\cdot}$ is the Kronecker delta function.  Then,
\begin{align}
\nonumber
\mC_\vtheta &= \bigl ( C_\vtheta(\vx_i, \vx_j) \bigr)_{i, j = 1}^n \\
\nonumber
& = \Bigl ( K_\vtheta \bigl(\vh(\phi(i-1) - \phi(j-1) \bigr) \bmod \vone ) \Bigr)_{i, j = 1}^n 
\qquad  \text{by \eqref{eqn:lattice_def} and \eqref{eq:shInv}}  \\
\nonumber
& = \biggl( 
\sum_{i',j'=1}^n \delta_{n\phi(i-1), i'-1}  \,
K_\vtheta \bigl (\vh (i'-j')/n \bmod \vone \bigr )
\delta_{j'-1,n\phi(j-1)} 
\biggr)_{i,j=1}^n \\
& = \mP \mK_{\vtheta}  \mP^T,  \label{shInvKernGramMatrix} 
\intertext{where } 
\mK_{\vtheta} &= \bigl ( K_\vtheta \bigl (\vh (i-j)/n \bmod \vone \bigr ) \bigr)_{i, j = 1}^n.
\end{align}

Because $\mK_\vtheta$ is circulant, we know the form of it's eigenvector-eigenvalue decomposition:
\begin{equation} \label{Keig}
    \mK_{\vtheta} = \frac 1n \mW \mLambda_\vtheta \mW^H, \; \ \mW =  \Bigl ( \me^{2 \pi \sqrt{-1} (i-1)(j-1)/n} \Bigr)_{i,j = 1}^n.
\end{equation}
By \eqref{shInvKernGramMatrix} we then have the eigenvector-eigenvalue decomposition for $\mC_{\vtheta}$ assumed in \eqref{eqn:ftk_factor}, namely
\begin{equation} \label{Clateig}
    \mC_{\vtheta} = \frac 1n \mV \mLambda_\vtheta \mV^H , \qquad \mV = \mP \mW \mP^T,
\end{equation}

where the eigenvalues of  $\mC_{\vtheta}$ and $\mK_{\vtheta}$ are identical.

Fast Bayesian transform assumptions \eqref{fastcompAssumpA} and \eqref{fastcompAssumpB} can be verified by \eqref{PermMat}, \eqref{Keig}, and \eqref{Clateig}.  Assumption \eqref{fastcompAssumpC} is satisfied because $\mV^H \vb = \mW^H \mP^T \vb$ is just the discrete Fourier transform of a vector whose rows have been permuted.  This can be performed in $\Order(n \log(n))$ operations by the FFT. A proper scaling of the kernel $K_\vtheta$, such as the one given by \eqref{the_kernel_eqn_bernoulli},  ensures that assumption \eqref{addAssump} is satisfied.

\subsection{Overcoming Cancellation Error}
For the covariance kernels used in our computation, it may happen that $n/\lambda_1$ is close to $1$.  Thus, the term $1-n/\lambda_1$, which appears in the credible interval half-widths, $\err_{\MLE}$, $\err_{\textup{full}}$, and $\err_{\textup{GCV}}$, may suffer from cancellation error and even become negative.  We have observed this phenomenon.  We can avoid this cancellation error by modifying how we compute the Gram matrix and its eigenvalues.

Any shift-invariant covariance kernel satisfying \eqref{addAssump} can be written as $C_\vtheta = 1 + \rC_\vtheta$, where $\rC_\vtheta$ is also symmetric and positive definite. The associated Gram matrix for $\rC_\vtheta$ is then $\rmC_\vtheta = \mC_\vtheta - \vone \vone^T$, and the eigenvalues of $\rmC_\vtheta$ are $\rlambda_1 = \lambda_1 - n, \lambda_2, \ldots, \lambda_n$, which follows because $\vone$ is the first eigenvector of both $\mC_\vtheta$ and $\rmC_\vtheta$. Then, 
\begin{equation*}
1 - \frac{n}{\lambda_1}  = \frac{\lambda_1 - n}{\lambda_1} = \frac{\rlambda_1}{\rlambda_1 +n},
\end{equation*}
where now the right hand side is free of cancellation error.

The covariance kernels that we use are of product form, namely,
\begin{equation*}
C_\vtheta(\vt, \vx) = \prod_{\ell=1}^d \left[1 + \rC_{\vtheta,\ell}(t_\ell,x_\ell) \right], \qquad  \rC_{\vtheta,\ell}:[0,1]^2 \to \reals.
\end{equation*}
Direct computation of $\rC_\vtheta(\vt,\vx) = C_\vtheta(\vt,\vx) -1$ introduces cancellation error if the $ \rC_{\vtheta, \ell}$ are small.  So, we employ the iteration
\begin{align*}
\rC^{(1)}_\vtheta(\vt,\vx) &= \rC_{\vtheta, 1}(t_1,x_1),  \\
\rC^{(\ell)}_\vtheta(\vt,\vx) &  = \rC^{(\ell-1)}_\vtheta(\vt,\vx) [1 + \rC_{\vtheta,\ell}(t_\ell,x_\ell)] + \rC_{\vtheta,\ell}(t_\ell,x_\ell), \\
& \hspace{5cm} \ell = 2, \ldots, d, \\
\rC_\vtheta(\vt,\vx)  & = \rC^{(d)}_\vtheta(\vt,\vx).
\end{align*}
In this way, the Gram matrix $\rmC_\vtheta$, whose $i,j$-element is $\rC_\vtheta(\vx_i,\vx_j)$ can be constructed in a way that avoids significant cancellation error.

Computing the eigenvalues of $\rmC$ via the procedure given in \eqref{eqn:fast_transform_to_eigvalues} yields $\rlambda_1 = \lambda_1 - n, \lambda_2, \ldots, \lambda_n$. The widths of the credible intervals in \eqref{errSimple} become
\begin{subequations}
\label{fastStoppingCriterions}
	\begin{align}
\label{fastStoppingCriterionMLE}
\err_\MLE  &
=
\frac{2.58}{n}\sqrt{
	\frac{\rlambda_1}{\lambda_1}
	\sum_{i=2}^{n} \frac{\abs{\widetilde{y}_i}^2}{\lambda_i}  
}, 
\\
\label{fastStoppingCriterionFull}
\err_{\textup{full}} 
& = \frac{t_{n-1,0.995}}{n} \sqrt{
	\frac{\rlambda_1}{n-1} \sum_{i=2}^n \frac{\abs{\widetilde{y}_i}^2}{\lambda_i}
}, \\
\label{fastStoppingCriterionGCV}
\err_{\textup{GCV}} & =
\frac{2.58}{n}\sqrt{	\frac{\rlambda_1}{\lambda_1} \sum_{i=2}^n \frac{\abs{\widetilde{y}_i}^2}{\lambda_i^2}  \left [ \frac 1n \sum_{i=1}^n \frac{1}{\lambda_i} \right]^{-1}} .
	\end{align}
\end{subequations}
For large $n$, $\lambda_1 \sim n$ and it follows that $\rlambda_1/\lambda_1 \approx \rlambda_1/(n-1)$  is small.  Moreover, for large $n$, the credible intervals via empirical Bayes and full Bayes are similar, since $t_{n-1,0.995} \approx 2.58$. 
The computational steps for the improved, faster, automatic Bayesian cubature are detailed in Algorithm \ref{algorithm2}.

\begin{algorithm}
	\caption{Fast Automatic Bayesian Cubature}\label{algorithm2}
	\begin{algorithmic}[1]
		\Require a generator for the rank-1 Lattice sequence
		$\vx_1, \vx_2, \ldots$; 
		a shift-invariant  kernel, $C_\vtheta$;
		a black-box function, $f$; 
		an absolute error tolerance,
		$\varepsilon>0$; the positive initial sample size, $n_0$, that is a power of $2$;
		the maximum sample size $n_{\textup{max}}$
		
		\State $n \gets n_0, \; n' \gets 0, \; \err_{\textup{CI}} \gets \infty$
		
		\While{$\err_{\textup{CI}} > \varepsilon$ and $n \le n_{\textup{max}}$}
		
		\State Generate $\{ \vx_i\}_{i=n' + 1}^{n}$ and sample $\{f(\vx_i)\}_{i=n'+1}^{n}$
		\State Compute $\vtheta$ by \eqref{eqn_MLE_loss_func_optimized_2} or \eqref{thetaGCV}
		\State Compute $\err_{\textup{CI}}$  according to \eqref{fastStoppingCriterionMLE}, \eqref{fastStoppingCriterionFull}, or \eqref{fastStoppingCriterionGCV}
		
		\State	$n' \gets n, \; n \gets 2n'$
		
		\EndWhile
		
		\State Update sample size to compute $\hmu$, $n \gets n'$
		\State Compute $\hmu$, the approximate integral,   according to \eqref{muhatGCV-FB-MLE-Simple}
		\State \Return $\hmu, \; n$  and $\err_{\textup{CI}}$
	\end{algorithmic}
\end{algorithm}

We summarize the results of this section and the previous one in the theorem below.  In comparison to Algorithm \ref{algorithm1}, the second and third components of the computational cost of Algorithm  \ref{algorithm2} are substantially reduced.
\begin{theorem}
Let $C_\vtheta$ be any symmetric, positive definite,  shift-invariant covariance kernel of the form \eqref{eq:shInv}, where $K_\vtheta$ has period one in every variable.  Furthermore, let $K_\vtheta$ be scaled to satisfy \eqref{addAssump}.  When matched with rank-1 lattice data sites, $C_\vtheta$ must satisfy fast Bayesian transform assumptions \eqref{fastcompAssump}.  The cubature, $\hmu$, is just the sample mean.  The \emph{fast Fourier transform} (FFT) can be used to expedite the estimates of $\vtheta$ in \eqref{fastTheta} and the credible interval half-widths \eqref{fastStoppingCriterions} so that Algorithm \ref{algorithm2} has a computational cost which is the sum of the following:
\begin{itemize}
	\item $\Order\bigl(n\$(f) \bigr)$ for the integrand data, where $\$(f)$ is the computational cost of a single $f(\vx)$;
	
	\item $\Order\bigl(N_{\opt} n \$(C_\vtheta) \bigr)$ for the evaluations of the vector $\vC_{1}$, where $N_\opt$ is the number of optimization steps required, and  $\$(C_\vtheta)$ is the computational cost of a single $C_\vtheta(\vt,\vx)$; and
	
	\item $\Order\bigl(N_{\opt} n \log(n) \bigr)$ for the FFT calculations; there is no $d$ dependence in these calculations.
	
\end{itemize}

\end{theorem}

Although the third part of the computational cost has the largest dependence on $n$, in practice it need not be the largest contributor to the computational cost.  If function values are the result of an expensive simulation, then the first part may consume most of the computation time.

We have implemented the fast adaptive Bayesian cubature algorithm in MATLAB as part of the Guaranteed Adaptive Integration Library (GAIL) \citep{ChoEtal19a} as \code{cubBayesLattice\_g}. This algorithm uses the covariance kernel defined in \eqref{eqn:ourkernel} with  $r=1$ and $2$, and the periodizing variable transforms in  \secref{period_var_tx}.  The rank-1 lattice node generator is taken from \cite{Nuy17a} \\ (\code{exod2\_base2\_m20}).

\section{Numerical Experiments} \label{sec:NumExp}

\subsection{Periodizing Variable Transformations}
\label{period_var_tx}
The shift-invariant covariance kernels underlying our Bayesian cubature  assume that the integrand has a degree of periodicity, with the smoothness assumed depending on the smoothness of the covariance kernel.  While integrands arising in practice may be smooth, they might not be periodic.  Variable transformations can be used to ensure periodicity.

Suppose that the original integral has been expressed as 
\begin{equation*}
\mu := \int_{[0,1]^d} g(\vt) \, \dif \vt,
\end{equation*}
where $g$ has sufficient smoothness, but lacks periodicity.  
The goal is to transform the integral above to the form of \eqref{eqn:defn_mu}, where the integrand $f$---and perhaps its derivatives---are  periodic.  

The baker's transform, also called the tent transform,
\begin{align} \label{eq:bakerTrans}
\vPsi: \vx \mapsto (\Psi(x_1),  \ldots, \Psi(x_d)),  \Psi(x)  =1 - 2 \abs{x - 1/2},
\end{align}
allows us to write $\mu$ in the form of \eqref{eqn:defn_mu}, where $f(\vx) = g(\vPsi(\vx))$.  Since  $\Psi'(x)$ is not continuous, $f$ does not have continuous derivatives. 

A family of variable transforms that can also preserve continuity of the derivatives from the original integrand $g$ takes the form
\begin{subequations} 
\begin{equation*}
\vPsi: \vx \mapsto (\Psi(x_1),  \ldots, \Psi(x_d)), \quad \Psi:[0,1] \mapsto [0,1].
\end{equation*}
This allows us to write $\mu$ in the form of \eqref{eqn:defn_mu} with
\begin{equation*}
f(\vx) = g(\vPsi(\vx)) \prod_{\ell = 1}^d \Psi'(x_l).
\end{equation*}
\end{subequations}
For $r \in \natzero$, if the following hold:
\begin{itemize}
	\item $\Psi \in C^{r+1}[0,1]$,
	\item  $\lim_{x \downarrow 0}x^{-r-1}\Psi'(x) = \lim_{x \uparrow 1} (1-x)^{-r-1}\Psi'(x) = 0$, and 
	\item $g \in C^{(r, \ldots, r)}[0,1]^d$,
\end{itemize}
then $f$ has continuous, periodic mixed partial derivatives of up to order $r$ in each direction.  
Examples of this kind of transform include \citep{Sid08a}:
\begin{align*}
\text{Sidi's } C^1 & : \Psi(x) = x - \frac{\sin(2\pi x)}{2 \pi},  \qquad   \Psi'(x) = 1 - \cos(2\pi x), \\
\text{Sidi's } C^2 & : \Psi(x) = \frac {8 - 9 \cos(\pi x) + \cos(3 \pi x)}{16} ,  
 \quad \Psi'(x) = \frac {3 \pi[3 \sin(\pi x) - \sin(3 \pi x)]}{16}.
\end{align*}

Periodizing variable transforms are used in the numerical examples below. In some cases, they can speed the convergence of the Bayesian cubature because they allow one to take advantage of smoother covariance kernels. 
However, there is a trade-off.  Smoother periodizing transformations tend to give integrands $f$ with larger inferred $s$ values and thus wider credible intervals.

\subsection{Test Results and Observations}

\begin{figure}
	\centering
	\includegraphics[width=0.98\linewidth]{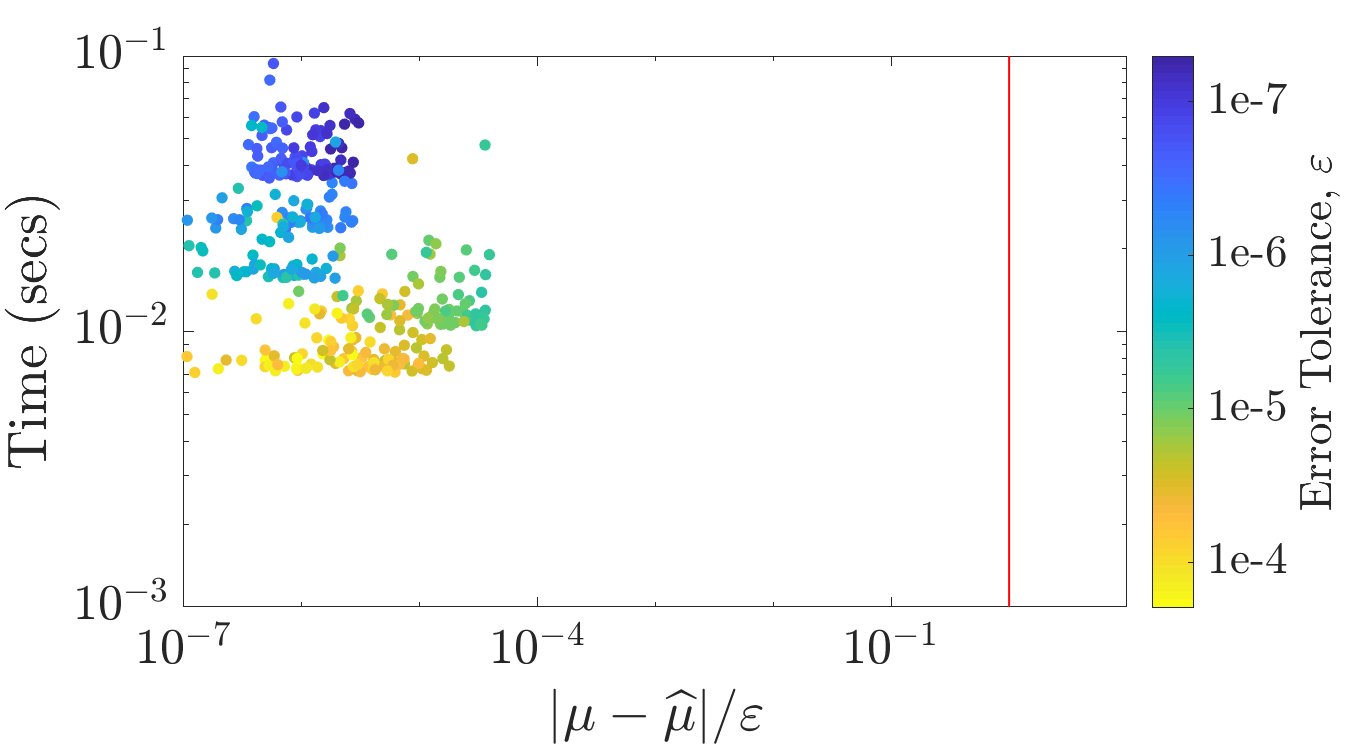}
	\caption[Guaranteed:]{Multivariate Gaussian probability example using the empirical Bayes stopping criterion.}
	\label{fig:mvn-guaranteed-MLE}
	\centering
	\includegraphics[width=0.98\linewidth]{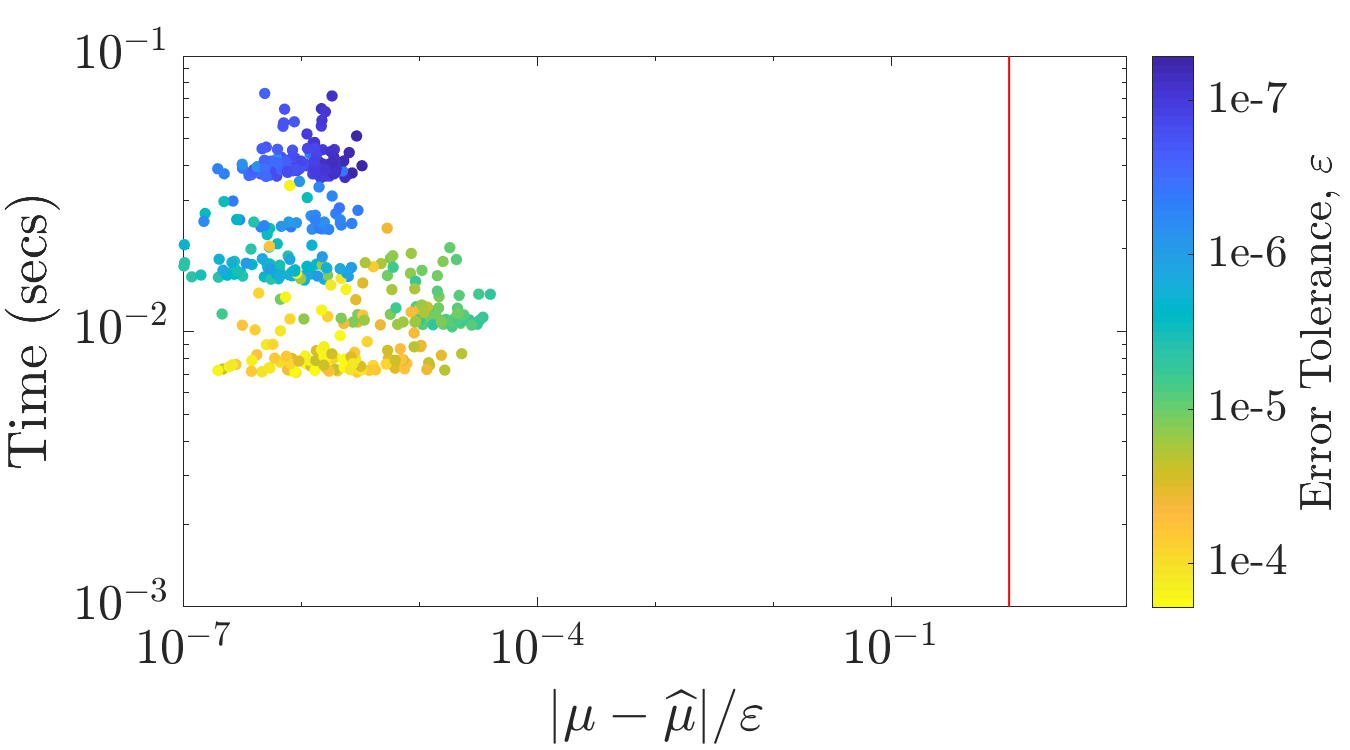}
	\caption[MVN guaranteed : FB]{Multivariate Gaussian probability example using the full Bayes stopping criterion.}
	\label{fig:mvn-guaranteed-FB}
	\centering
	\includegraphics[width=0.98\linewidth]{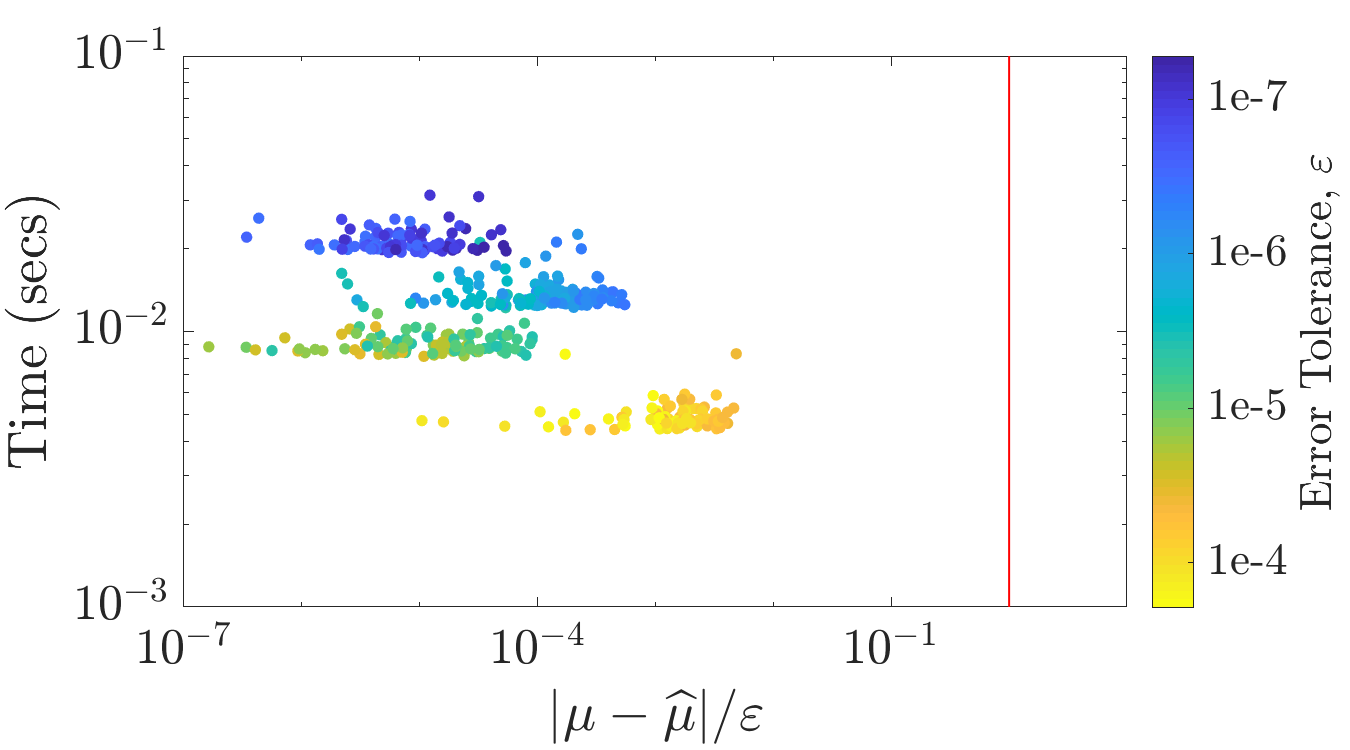}
	\caption[MVN guaranteed : GCV]{Multivariate Gaussian probability example using the GCV stopping criterion.}
	\label{fig:mvn-guaranteed-GCV}
\end{figure}

\begin{figure}
	\centering
	\includegraphics[width=0.95\linewidth]{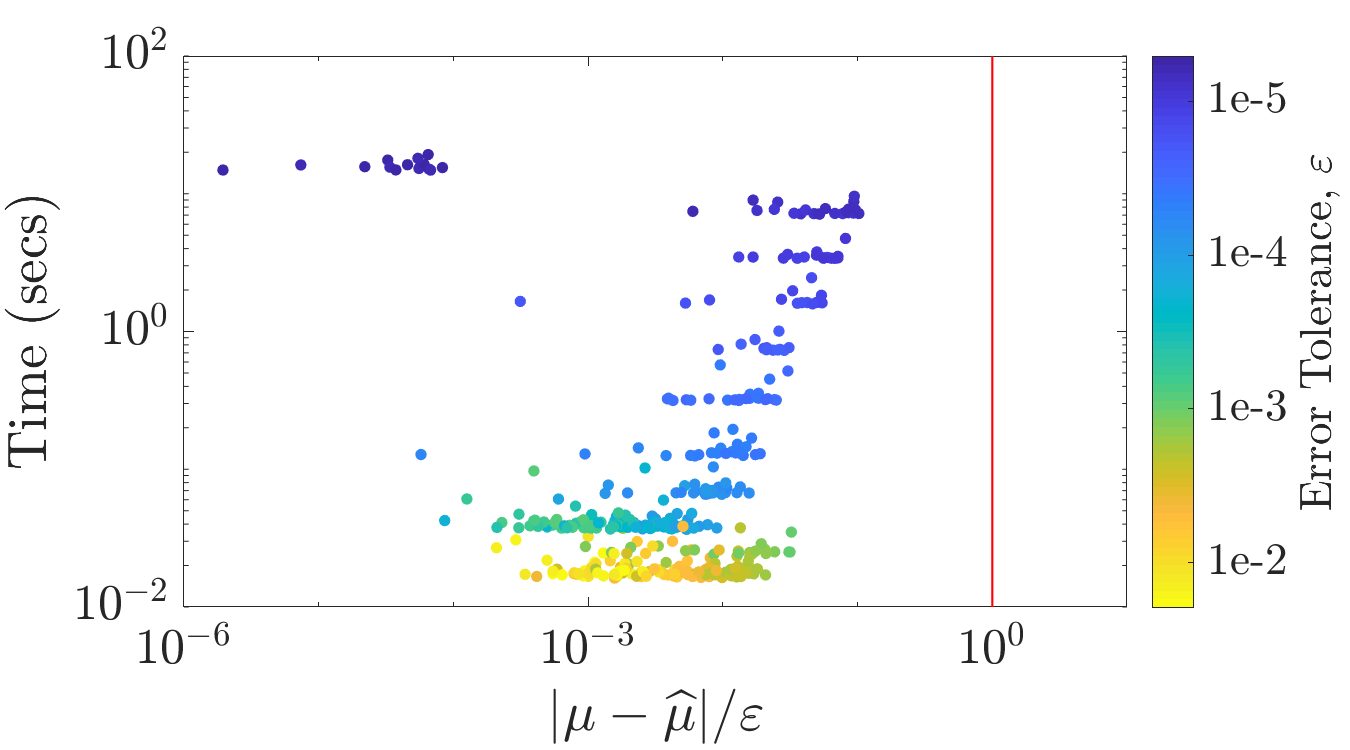}
	\caption[Keister guaranteed:MLE]{Keister's example using the empirical Bayes stopping criterion.}
	\label{fig:keister-guaranteed-MLE}
	\centering
	\includegraphics[width=0.95\linewidth]{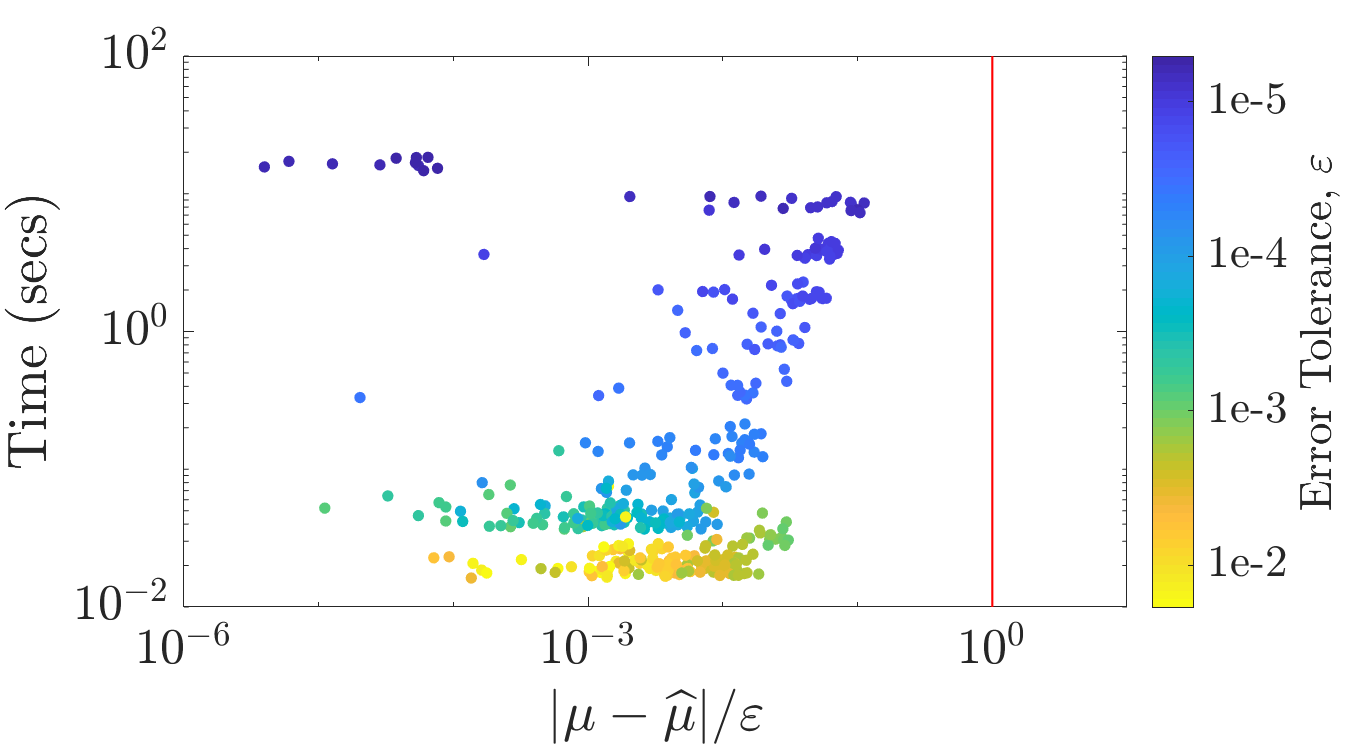}
	\caption[Keister guaranteed:FB]{Keister's example using the full Bayes stopping criterion.}
	\label{fig:keister-guaranteed-FB}
	\centering
	\includegraphics[width=0.95\linewidth]{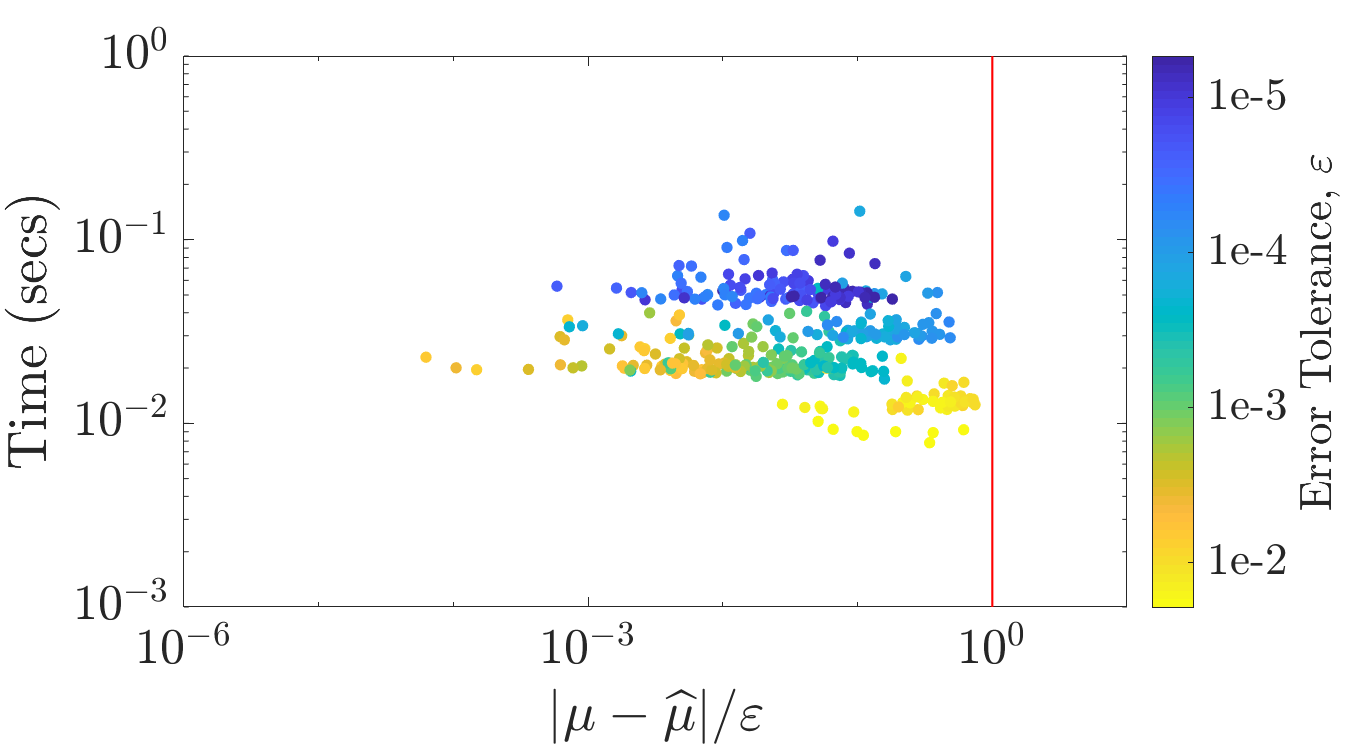}
	\caption[Keister guaranteed:GCV]{Keister's example using the GCV stopping criterion.}
	\label{fig:keister-guaranteed-GCV}
\end{figure}

\begin{figure}
	\centering
	\includegraphics[width=0.95\linewidth]{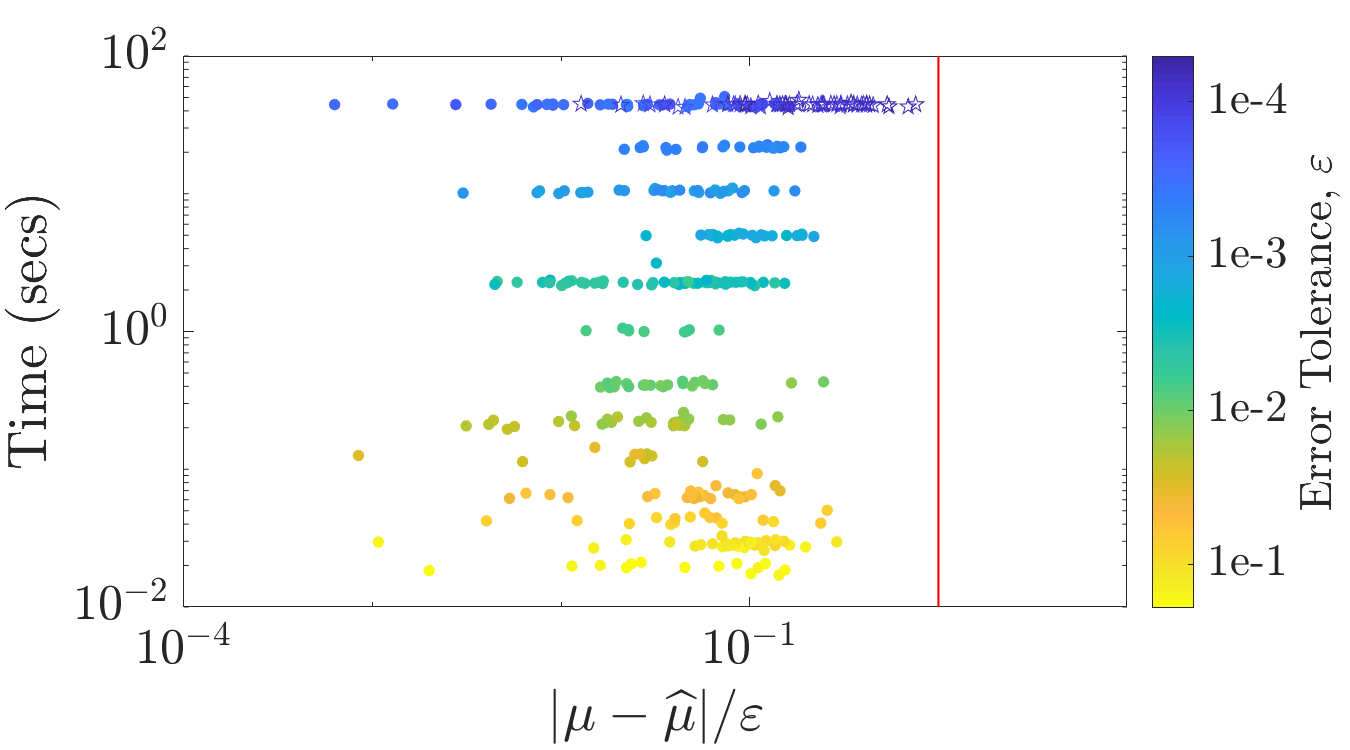}
	\caption[Option pricing Guaranteed: MLE]{Option pricing using the empirical Bayes stopping criterion.}
	\label{fig:optprice-guaranteed-MLE}
	\centering
	\includegraphics[width=0.95\linewidth]{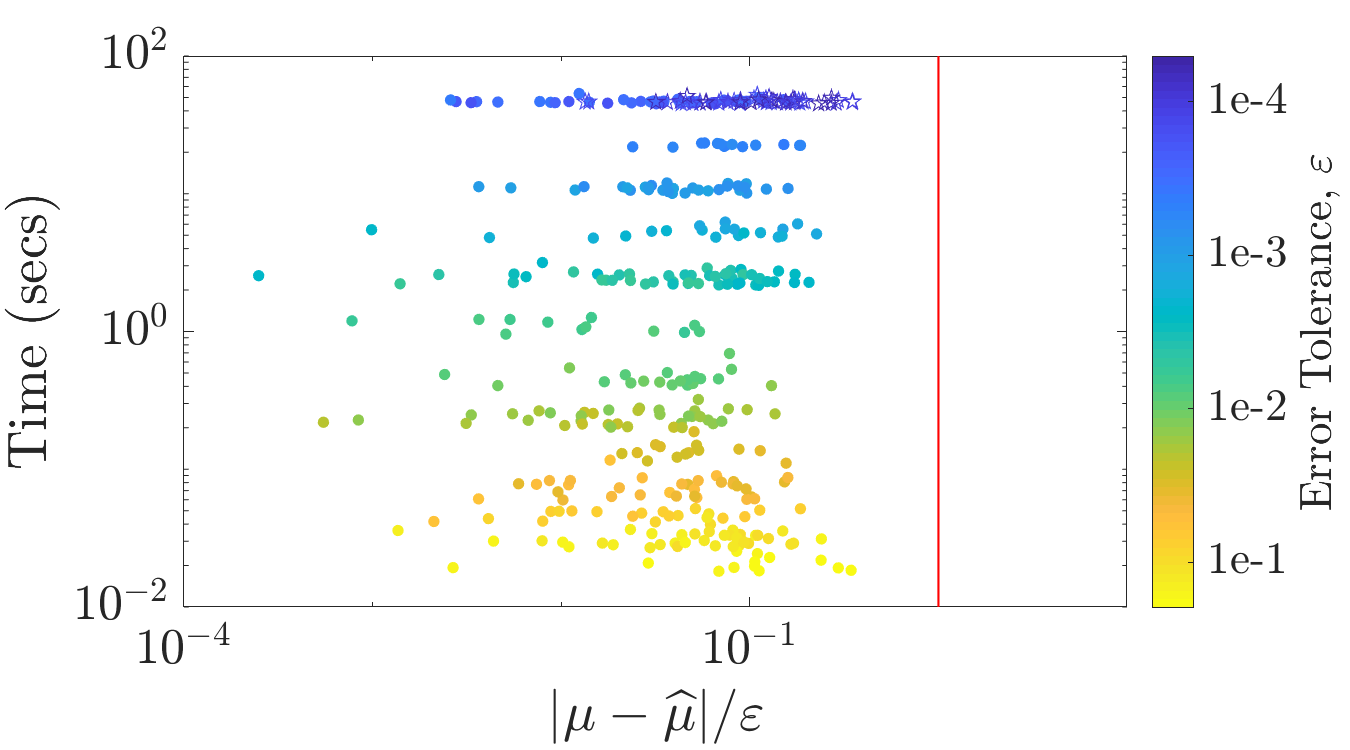}
	\caption[OptPrice guaranteed : FB]{Option pricing using the full Bayes stopping criterion.}
	\label{fig:optprice-guaranteed-FB}
	\centering
	\includegraphics[width=0.95\linewidth]{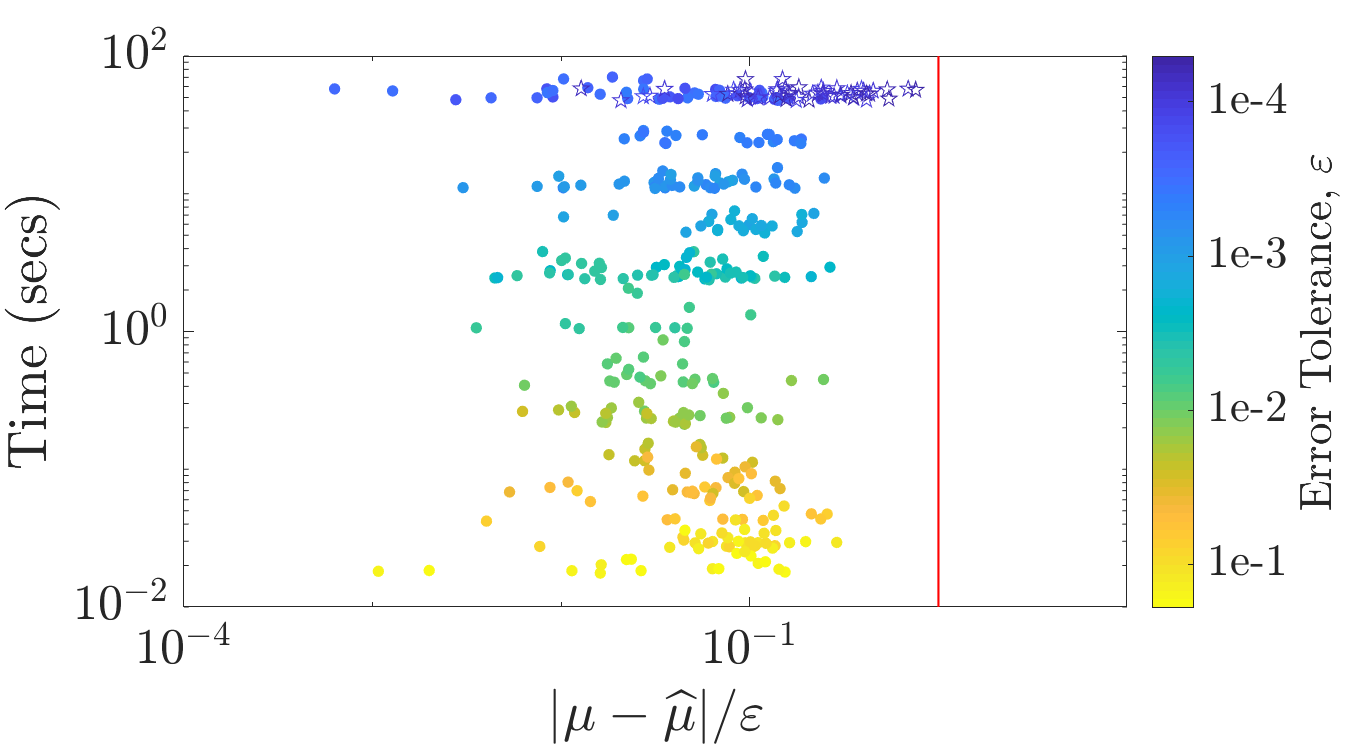}
	\caption[OptPrice guaranteed : GCV]{Option pricing using the  GCV stopping criterion.}
	\label{fig:optprice-guaranteed-GCV}
\end{figure}

Three integrals were evaluated using the GAIL algorithm \code{cubBayesLattice\_g}:  a multivariate Gaussian  probability, Keister's example, and an option pricing example. 
These three integrands are defined below.
The sequences $\{\vx_i\}_{i=1}^\infty$ are the randomly shifted lattice node sequences supplied by GAIL. 
For each integral and each of our stopping criteria---empirical Bayes, full Bayes, and generalized cross-validation---our algorithm was run for $400$ different randomly chosen error tolerances. The error tolerances were chosen randomly in an interval depending on the difficulty of the problem. In each run, the nodes were also randomly shifted with $\mathcal{U}[0,1]$ shifts independent of each other and the error tolerances. The accuracy of the algorithm depends mildly on the shift; there is no universally optimal shift.
For each test, the execution times are plotted against $\abs{\mu - \hmu}/\varepsilon$.  We expect $\abs{\mu - \hmu}/\varepsilon$ to be no greater than one, but hopefully not too much smaller than one, which would indicate a stopping criterion that is too conservative. Figures \ref{fig:mvn-guaranteed-MLE} to \ref{fig:optprice-guaranteed-GCV} can be reproduced using the script \code{cubBayesLattice\_guaranteed\_plots.m} in GAIL.

Ideally, we would optimize both $r$ and $\eta$ simultaneously in the definition of our $C_\vtheta$ in \eqref{eqn:ourkernel}.  However, in these examples we fix $r$ and optimize $\eta$ only.  This is a technical challenge, not a limitation of our theory.

\paragraph{Multivariate Gaussian Probability.}

This example was already introduced in Section \ref{MVN_example}, where we used the Mat\'ern covariance kernel.  Here we apply Sidi's $C^2$  periodization to $ f_{\textup{Genz}}$ \eqref{fGenzdef} and choose $d'=3$, $d=2$, and $r=2$. The simulation results for this example are summarized in Figures \ref{fig:mvn-guaranteed-MLE}, \ref{fig:mvn-guaranteed-FB}, and \ref{fig:mvn-guaranteed-GCV}.  In all cases, the Bayesian cubature returns an approximation within the prescribed error tolerance. We used the same setting as before with generic slow Bayesian cubature in \secref{MVN_example} for comparison. For error tolerance $\varepsilon=10^{-5}$ with the empirical Bayes stopping criterion, our fast algorithm takes just under 0.01 seconds as shown in \figref{fig:mvn-guaranteed-MLE} whereas the generic algorithm takes over 20 seconds as shown in \figref{fig:MVN_Metern_d2b2}. 

Amongst the three stopping criteria, GCV achieves the desired tolerance faster than the others. 
One can also observe from the figures, the credible intervals are in general much wider than the true error.
This could be due to the periodized integrand being smoother than the $r=2$ covariance kernel assumes. Perhaps one should consider smoother covariance kernels.

\paragraph{Keister's Example.}

This multidimensional integral function comes from \cite{Kei96} and is inspired by a physics application:
\begin{align*}
\mu & =  \int_{\reals^d} \cos(\norm{ \vt}) \exp(-\norm{ \vt }^2) \, \dvt 
 = \int_{[0,1]^d} f_{\textup{Keister}}(\vx) \, \dvx,\\
\text{where } 
f_\textup{Keister}(\vx) &= \pi^{d/2} \cos\left(\norm{ \Phi^{-1}(\vx)/2}\right)  ,
\end{align*}
and again $\Phi$ is the standard Gaussian distribution.
The true value of $\mu$ can be calculated iteratively in terms of a quadrature as follows:  
\begin{equation*}
\mu = \frac{2 \pi^{d/2} I_c(d)}{\Gamma(d/2)}, \quad d=1,2, \ldots
\end{equation*}
where $\Gamma$ denotes the gamma function, and
\begin{align*}
I_c(1) &= \frac{\sqrt{\pi}}{2 \exp(1/4)}, 
\\
I_s(1) &= \int_{x=0}^\infty \exp(-\vx^T\vx)\sin(\vx) \, \dvx 
\; = \;  0.4244363835020225,
\\
I_c(2) &= \frac{1-I_s(1)}{2}, \qquad
I_s(2) = \frac{I_c(1)}{2}
\\
I_c(j) &= \frac{(j-2)I_c(j-2)-I_s(j-1)}{2},
\qquad j =3,4,\ldots
\\
I_s(j) &= \frac{(j-2)I_s(j-2)-I_c(j-1)}{2},
\qquad j =3,4,\ldots.
\end{align*}

Figures \ref{fig:keister-guaranteed-MLE}, \ref{fig:keister-guaranteed-FB}, and \ref{fig:keister-guaranteed-GCV} summarize the numerical tests for this integral.  We used Sidi's $C^1$ periodization, dimension $d=4$, and $r=2$. 
As we can see, the GCV stopping criterion achieves results faster than the other stopping criteria, similar to the multivariate Gaussian case. The credible intervals are narrower than the multivariate Gaussian case since the covariance kernel is smoother than the periodization transform used.

\paragraph{Option Pricing.}

The price of financial derivatives can often be modeled by high dimensional integrals. If the underlying asset is described in terms of a discretized geometric Brownian motion, then the fair price of the option is:
\begin{equation*}
\mu = \int_{\reals^d} \code{payoff}(\vz) \frac{\exp(\frac 12 \vz^T\mSigma^{-1}\vz)}{\sqrt{(2\pi)^d \det(\mSigma)}} \, \dvz = \int_{[0,1]^d} f(\vx) \, \dvx,
\end{equation*} 
where \code{payoff($\cdot$)} defines the discounted payoff of the option,
\begin{align*}
\mSigma &= (T/d) \bigl(\min(j,k) \bigr)_{j,k=1}^d = \mL \mL^T,\\
f(\vx) &= \code{payoff} \left(\mL 
\begin{pmatrix}
\Phi^{-1}(x_1) \\ \vdots \\ \Phi^{-1}(x_d)
\end{pmatrix} \right).
\end{align*}
The Asian arithmetic mean call option has a payoff of the form
\begin{align*}
\code{payoff}(\vz) &= \max\left( \frac 1d  \sum_{j=1}^d S_j(\vz) - K, 0 \right) \me^{-R T}, \\
\text{where}\;
S_j(\vz) &= S_0 \exp\bigl((R-\sigma^2/2)j(T/d) + \sigma \sqrt{(T/d)} z_j \bigr).
\end{align*}
Here, $T$ denotes the time to maturity of the option, $d$ the number of time steps, $S_0$ the initial price of the stock, $R$ the interest rate, $\sigma$ the volatility, and $K$ the strike price.  

Figures \ref{fig:optprice-guaranteed-MLE}, \ref{fig:optprice-guaranteed-FB}, and 
\ref{fig:optprice-guaranteed-GCV} summarize the numerical results for this example using
$
T = 1/4, \ \ d = 13, \ \ S_0 = 100, \ \ R =  0.05, \ \ \sigma = 0.5, \ \ K = 100.
$
Moreover, $\mL$ is chosen to be the matrix of eigenvectors of $\mSigma$ times the square root of the diagonal matrix of eigenvalues of $\mSigma$. 
Because the integrand has a kink caused by the $\max$ function, it does not help to use a periodizing transform that is very smooth.  We choose the baker's transform \eqref{eq:bakerTrans} and $r = 1$.

In summary, the Bayesian cubature algorithm computes the integral within the user-specified tolerance in nearly all of the test cases.  The rare exceptions occurred in the option pricing example for $\varepsilon \le 10^{-4}$. Our algorithm used the maximum allowed sample size and still did not reach the stopping criterion $\err_{\textup{CI}} \leq \varepsilon$, due to the complexity and high dimension of the integrand. Those cases are shown as hollow stars in Figures \ref{fig:optprice-guaranteed-MLE}, \ref{fig:optprice-guaranteed-FB}, and 
\ref{fig:optprice-guaranteed-GCV}.

One may question whether an integrand with non-negative values is well represented by a Gaussian process.  Since we allow a nonzero mean, this assumption is somewhat more palatable. Bayesian algorithms assuming non-Gaussian processes are more difficult to execute, but this is an area for further research.

A noticeable observation from the plots in all three examples is how the ratio of the true error to the error tolerance varies from nearly one all the way down to $10^{-7}$.  Since the credible interval half-widths are not much smaller than $\varepsilon$, this means that the credible intervals are quite conservative in many cases.  For option pricing example, this is less of an issue than for the multivariate Gaussian and Keister's examples. The reason that the credible intervals wildy overestimate the true error for the multivariate Gaussian and Keister's examples may be that these integrands are significantly smoother than the assumed covariance kernel.  This is a matter for further investigation.

\section{Discussion and Further Work}

We have developed a fast, automatic Bayesian cubature that estimates a multidimensional definite integral within a user defined error tolerance.  The stopping criteria arise from assuming the integrand to be an instance of a Gaussian process. There are three approaches:  empirical Bayes, full Bayes, and generalized cross-validation.  The computational cost of the automatic Bayesian cubature can be dramatically reduced if the covariance kernel matches the nodes.  One such match in practice is rank-1 lattice nodes and shift-invariant kernels.  The matrix-vector multiplications can be accomplished using the fast Fourier Transform.  The performance of our automatic Bayesian cubature is illustrated using three integration problems.  

Digital sequences and digital shift and/or scramble invariant kernels have the potential of being another match that satisfies the conditions in Section \ref{sec:fast_BC}.  The fast Bayesian transform would correspond to a fast Walsh transform.  


One should be able to adapt our Bayesian cubature to control variates, i.e., assuming  
\begin{equation*}
f = \mathcal{GP} \left( \beta_0 + \beta_1 g_1 + \cdots + \beta_p g_p, s^2 C_\vtheta \right),
\end{equation*}
for some choice of $g_1, \ldots, g_p$ whose integrals are known, and some parameters $\beta_0, \ldots, \beta_p$ in addition to $s$ and $C_\vtheta$.  The efficacy of this approach has not yet been explored.




\begin{acknowledgements}
This research was supported in part by the National Science Foundation grants DMS-1522687 and DMS-1638521 (SAMSI).
The authors would like to thank the organizers of the SAMSI-Lloyds-Turing Workshop on Probabilistic Numerical Methods, where a preliminary version of this work was discussed.  The authors also thank Chris Oates and Sou-Cheng Choi for valuable comments.
\end{acknowledgements}

\bibliographystyle{abbrvnat}
\bibliography{FJHown23,FJH23}

\begin{appendices}
\section{Details of the Full Bayes Posterior Density for $\mu$} \label{appendix:full_bayes}
To simplify, we drop the dependence of $c_{0\vtheta}$, $\vc_\vtheta$, and $\mC_\vtheta$ on $\vtheta$ in the notation below.
Starting from the Bayesian formula for the posterior density for $\mu$ at the beginning of Section \ref{sec:fullBayes} with the non-informative prior, it follows that 
\begin{align*}
\rho_{\mu|\vf}(z | \vy)
& \propto \int_{0}^\infty \int_{-\infty}^\infty 
\rho_{\mu | m, s^2, \vf}(z | \xi, \lambda , \vy) 
  \rho_{\vf | m, s^2} (\vy | \xi, \lambda) \, \vrho_{m,s^2}(\xi,\lambda)  \, \D \xi \D \lambda \\
& \propto \displaystyle \int_{0}^\infty  \frac{1}{\lambda^{(n+3)/2}} 
 \int_{-\infty}^\infty  \exp \biggl( -\frac{1}{2\lambda}\biggl\{
\frac{
	[z - \xi (1 - \vc^T \mC^{-1} \vone)  -  \vc^T \mC^{-1} \vy]^2}
{c_0  -\vc ^T \mC^{-1} \vc}  \\
& \qquad + (\vy - \xi \vone)^T \mC^{-1}(\vy - \xi \vone) \biggr \} \biggr) \, \D \xi \D \lambda \\
& \qquad \qquad
\text{by \eqref{eqn:fGaussDist}, \eqref{eqn:condInteg}} \; \text{and} \; \rho_{m,s^2}(\xi,\lambda) \propto 1/\lambda \\
& \propto \displaystyle \int_{0}^\infty  \frac{1}{\lambda^{(n+3)/2}} 
\int_{-\infty}^\infty  \exp\left( -\frac{\alpha \xi^2 -2 \beta \xi + \gamma}{2\lambda(c_0  -\vc ^T \mC^{-1} \vc)} \right) \, \D \xi \D \lambda,
\end{align*}
where
\begin{align*}
\alpha & = (1 - \vc^T \mC^{-1} \vone)^2 + \vone^T \mC^{-1} \vone (c_0  -\vc ^T \mC^{-1} \vc),\\
\beta & =(1 - \vc^T \mC^{-1} \vone)(z - \vc^T \mC^{-1} \vy )  
+ \vone^T \mC^{-1} \vy (c_0  -\vc ^T \mC^{-1} \vc),\\
\gamma &  = (z - \vc^T \mC^{-1} \vy )^2  + \vy^T \mC^{-1} \vy (c_0  -\vc ^T \mC^{-1} \vc).
\end{align*}
In the derivation above and below, factors that are \emph{independent} of $\xi$, $\lambda$, or $z$ can be discarded since we only need to preserve the proportion.  But, factors that depend on $\xi$, $\lambda$, or $z$ must be kept.  
Completing the square, $
\alpha \xi^2 -2 \beta \xi + \gamma 
= \alpha (\xi -\beta/\alpha)^2  - (\beta^2/\alpha) + \gamma,
$
allows us to evaluate the integrals with respect to $\xi$ and $\lambda$:
\begin{align*}
\rho_{\mu|\vf}(z | \vy)
& \propto \displaystyle \int_{0}^\infty  \frac{1}{\lambda^{(n+3)/2}}  \exp\left( -\frac{  \gamma - \beta^2/\alpha}{2\lambda(c_0  -\vc ^T \mC^{-1} \vc)} \right)  \\
&\qquad \qquad \times
\int_{-\infty}^\infty  \exp\left( -\frac{\alpha (\xi -\beta/\alpha)^2}{2\lambda(c_0  -\vc ^T \mC^{-1} \vc)} \right) \, \D \xi \D \lambda \\
& \propto \displaystyle \int_{0}^\infty  \frac{1}{\lambda^{(n+2)/2}}  \exp\left( -\frac{  \gamma - \beta^2/\alpha}{2\lambda(c_0  -\vc ^T \mC^{-1} \vc)} \right) \D \lambda \\
& \propto \left(\gamma - \frac{\beta^2}{\alpha}\right)^{-n/2} \\
& \propto \left(\alpha \gamma - \beta^2\right)^{-n/2}.
\end{align*}
Finally, we simplify the key term via straightforward calculations to the following:
\begin{align*}
\alpha \gamma - \beta^2 \propto 1 +  \frac{(z - \hmu_{\MLE})^2}{(n-1)s_{\textup{full}}^2},
\end{align*}
where 
\begin{multline*}
\hsigma_{\textup{full}}^2
:= \frac{1}{n-1}
\vy^T\left[ \mC^{-1} 
- \frac{ \mC^{-1} \vone\vone^T \mC^{-1}}{\vone^T \mC^{-1} \vone}  \right]\vy
\left[\frac{(1 - \vc^T \mC^{-1} \vone)^2}{\vone^T \mC^{-1} \vone} + (c_0  -\vc ^T \mC^{-1} \vc) \right].
\end{multline*}
This completes the derivation of \eqref{eqn:sigma2_full}.
\end{appendices}

\end{document}